\documentclass[12pt,a4paper,reqno]{amsart}
\usepackage{amssymb}
\usepackage{latexsym}
\usepackage{amsmath}
\usepackage{mathrsfs}
\usepackage{cite}
\usepackage{textcomp}
\usepackage{upgreek}
\usepackage{mathabx}
\usepackage{amscd}
\usepackage{color}

\usepackage{calc}

\newcommand{\scal}[2]{\langle #1,#2\rangle}
\newcommand{\nn}[1]{\mathbf N^{#1}}
\newcommand{\rr}[1]{\mathbf R^{#1}}

\newcommand{\cc}[1]{\mathbf C^{#1}}
\newcommand{\zz}[1]{\mathbf Z^{#1}}
\newcommand{\nm}[2]{\Vert #1\Vert _{#2}}
\newcommand{\Nm}[2]{\left \Vert #1 \right \Vert _{#2}}

\newcommand{\op}{\operatorname{Op}}
\newcommand{\sets}[2]{\{ \, #1\, ;\, #2\, \} }
\newcommand{\Sets}[2]{\left \{ \, #1\, ;\, #2\, \right \} }
\newcommand{\ep}{\varepsilon}

\newcommand{\cdo}{\, \cdot \, }

\newcommand{\eabs}[1]{\langle #1\rangle}     

\newcommand{\vrum}{\vspace{0.2cm}}

\newcommand{\maclA}{\mathcal A}

\newcommand{\maclH}{\mathcal H}

\newcommand{\maclS}{\mathcal S}

\newcommand{\mascF}{\mathscr F}
\newcommand{\mascP}{\mathscr P}
\newcommand{\mascS}{\mathscr S}

\newcommand{\mabfp}{{\boldsymbol p}}
\newcommand{\mabfq}{\boldsymbol q}
\newcommand{\mabfr}{\boldsymbol r}

\newcommand{\GL}{\mathbf{M}}

\numberwithin{equation}{section}          

\newtheorem{thm}{Theorem}
\numberwithin{thm}{section}
\newtheorem{prop}[thm]{Proposition}

\newtheorem{lemma}[thm]{Lemma}

\newcommand{\rubrik}{}

\theoremstyle{definition}

\newtheorem{defn}[thm]{Definition}

\theoremstyle{remark}
\newtheorem{rem}[thm]{Remark}

\newcommand{\rd}{\mathbf{R} ^{d}}

\title{Pseudo-differential calculus in a Bargmann setting}

\author{Nenad Teofanov}

\address{Department of Mathematics and Informatics,
University of Novi Sad, Novi Sad, Serbia}

\email{nenad.teofanov@dmi.uns.ac.rs}

\author{Joachim Toft}

\address{Department of Computer science, Mathematics and Physics,
Linn{\ae}us University, V{\"a}xj{\"o}, Sweden}

\email{joachim.toft@lnu.se}

\keywords{Analytic kernels, Berezin operators, Pilipovi{\'c} spaces,
modulation spaces, Gelfand-Shilov spaces}

\subjclass[2010]{Primary: 32W25, 35S05, 32A17, 46F05, 42B35
\quad Secondary: 32A25, 32A05}

\begin{document}

\begin{abstract}
We give a fundament for Berezin's analytic $\Psi$do
considered in \cite{Berezin71} in terms of Bargmann images
of Pilipovi{\'c} spaces. We deduce basic continuity results for such
$\Psi$do, especially when the operator kernels are in suitable
mixed weighted Lebesgue spaces and
act on certain weighted Lebesgue spaces of
entire functions. In particular, we show how these results imply well-known
continuity results for real $\Psi$do with symbols in modulation
spaces, when acting on other modulation spaces.
\end{abstract}

\maketitle

\section{Introduction}\label{sec0}

\par

The aim of the paper is to put a fundament for the theory of analytic
pseudo-differential operators, considered in \cite{Berezin71} by F. Berezin.
This is essentially done through a detailed analysis  of Bargmann images of the so-called
Pilipovi{\'c} spaces of functions and distributions, given in \cite{FeGaTo2,Toft18}.
More precisely, we consider kernels related to integral representations of
analytic pseudo-differential operators to deduce their continuity properties.
When the corresponding symbols belong to suitable (weighted) Lebesgue spaces of
semi-conjugate analytic functions, we prove the continuity of the analytic
pseudo-differential operators  when acting between (weighted) Lebesgue spaces
of analytic functions. Moreover, by using the relationship between the Bargmann
transform and the short-time Fourier transform we show that our results can be used
to recover well-known (sharp) continuity
properties of (real) pseudo-differential operators with symbols in modulation spaces
which act between other modulation spaces, see \cite{To8,To11,Toft22}.
We emphasize that our approach here is more general, because we have
relaxed the assumptions on the involved weight functions, compared to
earlier contributions.

\par

Analytic pseudo-differential operators, considered in \cite{Berezin71}
by Berezin are well-designed when considering several problems in
analysis and its applications, e.{\,}g. in quantum mechanics. In the
context of abstract harmonic analysis it follows that any linear and continuous
operator between Fourier invariant function and (ultra-)distribution spaces may,
in a unique way, be transformed into an analytic pseudo-differential operator by
the Bargmann transform (see Section \ref{sec2}). An advantage of
such reformulations is that all of the involved objects are essentially entire
functions and thereby possess several strong and convenient
properties.

\par

The definition of analytic pseudo-differential operators resembles the
definition of real pseudo-differential operators. In fact, let $a(x,\xi )$ be a
suitable function or (ultra-)distribution on the phase space $\rr {2d}$.
Then the (real) pseudo-differential operator $\op (a)$ acting on
suitable sets of functions or (ultra-)distributions
on the configuration space $\rr d$ is given by
\begin{equation}\label{Eq:RealPseudo}
f(x)\mapsto (\op (a)f)(x) = (2\pi )^{-\frac d2}\int _{\rr d}
a(x,\xi )\widehat f(\xi )e^{i\scal x\xi}\, d\xi .
\end{equation}
Here the integral in \eqref{Eq:RealPseudo} should
be interpreted in a distributional (weak) sense, if
necessary, and we refer to \cite{Ho1}
or Section \ref{sec1} for the notation.

\par

Suppose instead that $a$ is a suitable semi-conjugate entire (analytic)
function on $\cc d\times \cc d
\asymp \cc {2d}$, i.{\,}e. $(z,w)\mapsto a(z,\overline w)$ is an
entire (analytic) function.
Then the analytic pseudo-differential operator $\op _{\mathfrak V}(a)$
acting on suitable entire functions $F$ on $\cc d$ is given by
\begin{equation}\label{Eq:CompPseudo1}
F(z)\mapsto (\op _{\mathfrak V}(a)F)(z) = \int _{\cc d}
a(z,w)F(w)e^{(z,w)}\, d\mu (w).
\end{equation}
Here $d\mu (w)$ is the Gauss measure $\pi ^{-d}e^{-|w|^2}\, d\lambda (w)$, where
$d\lambda (w)$ is the Lebesgue measure on $\cc d$, and $ (z,w) = \sum_{j=1} ^{d} z_j \cdot \overline {w_j}, $ when
$z=(z_1,\dots ,z_d) \in \cc d$, and $ w=(w_1,\dots ,w_d)\in \cc d$. 
This means that the operator kernel (with respect to $d\mu$) is given by
\begin{equation}\label{Eq:CompPseudo2}
K(z,w)=K_a(z,w) = a(z,w)e^{(z,w)}.
\end{equation}
Evidently, $(\op _{\mathfrak V}(a)F)(z)$ is equal to the integral
operator
\begin{equation}\label{Eq:CompKernelOp}
(T_KF)(z) = \int _{\cc d}
K(z,w)F(w)\, d\mu (w)
\end{equation}
with respect to $d\mu$, when $K$ is given by \eqref{Eq:CompPseudo2}.
By the analyticity properties of the symbol $a$ it follows that
$(z,w)\mapsto K(z,\overline w)$ is an entire function on $\cc {2d}$.

\par

In \cite{Berezin71,To11} several facts
of analytic pseudo-differential operators are deduced. For example,
if $a$ and $F$ are chosen such that
$$
z\mapsto a(z,\cdo )F\cdot e^{(z,\cdo )}
$$
is locally uniformly bounded
and analytic from $\cc d$ to $L^1(d\mu )$, then
$\op _{\mathfrak V}(a)F$ in \eqref{Eq:CompPseudo1} is a well-defined
entire function on $\cc d$.
In \cite{Berezin71,To11} it is also observed that
\begin{equation}\label{Eq:BargmannCreationAnnihilation}
(\op _{\mathfrak V}(z_j)F)(z) = z_jF(z)
\quad \text{and}\quad
(\op _{\mathfrak V}(\overline w_j)F)(z) = (\partial _{j}F)(z)
\end{equation}
when $F \in L^1(d\mu _1)\cap A(\cc d)$, and $d\mu _1(w) = (1+|w|)\, d\mu (w)$.

\medspace

In such setting we study the mapping properties for complex integral operators and
pseudo-differential operators when respectively $K=K_a$ and
$a$ above belong to suitable classes of semi-conjugate entire functions.
In fact, we permit more generally that $K$ and $a$ belong to suitable classes
of formal semi-conjugate analytic power series expansions.
That is, $K(z,w)$ and $a(z,w)$ are of the forms
$$
\sum _{\alpha ,\beta}c_K (\alpha ,\beta) e_\alpha (z)e_{\beta}(\overline w)
\quad \text{and}\quad
\sum _{\alpha ,\beta}c_a (\alpha ,\beta) e_\alpha (z)e_{\beta}(\overline w),
\qquad e_\alpha (z)=\frac {z^\alpha}{\sqrt {\alpha !}},
$$
respectively.

\par

To set the stage for our study we collect the background material in Section \ref{sec1}.
It contains a brief account on weight functions,  Gelfand-Shilov spaces, spaces
of Hermite functions and power series expansions,
modulation spaces, and Bargmann transform and spaces of analytic functions.
Especially, we recall basic facts for the spaces
\begin{alignat}{2}
& \maclA _{\flat _\sigma}(\cc d) \  \big (\maclA _{0,\flat _\sigma}(\cc d) \big ),&
\quad
& \maclA _{s}(\cc d) \  \big ( \maclA _{0,s}(\cc d) \big ),
\label{Eq:AnalSpaces1}
\intertext{and their \emph{Bargmann} duals}
& \maclA _{\flat _\sigma}'(\cc d)  \ \big ( \maclA _{0,\flat _\sigma}'(\cc d)\big ),&
\quad
& \maclA _{s}'(\cc d) \ \big ( \maclA _{0,s}'(\cc d) \big ),
\label{Eq:AnalSpaces2}
\end{alignat}
when $s,\sigma >0$. The spaces in \eqref{Eq:AnalSpaces1} consist of all
formal power series
\begin{equation}\label{Eq:FormalPowerSeriesIntrod}
F(z) = \sum _{\alpha}c(F,\alpha) e_\alpha (z),
\end{equation}
with coefficients satisfying
\begin{alignat*}{2}
|c(F,\alpha )| &\lesssim h^{|\alpha |}\alpha !^{-\frac 1{2\sigma}},
\quad
|c(F,\alpha )| &\lesssim e^{-r|\alpha |^{\frac 1{2s}}},
\intertext{respectively, for some (for every) $h,r>0$, and the
spaces in \eqref{Eq:AnalSpaces2} consist of all
formal power series in \eqref{Eq:FormalPowerSeriesIntrod} such that}
|c(F,\alpha )| &\lesssim h^{|\alpha |}\alpha !^{+\frac 1{2\sigma}},
\quad
|c(F,\alpha )| &\lesssim e^{+r|\alpha |^{\frac 1{2s}}},
\end{alignat*}
respectively, for every (for some) $h,r>0$.

\par

In Section \ref{sec2} we extend the definition of \eqref{Eq:CompKernelOp}
to allow the kernels $K$ to belong to any of the spaces
\begin{alignat}{4}
&\wideparen \maclA _{0,\flat _\sigma}(\cc {2d}),
&\quad
&\wideparen \maclA _{\flat _\sigma}(\cc {2d}),
&\quad
&\wideparen \maclA _{0,s}(\cc {2d}),
&\quad
&\wideparen \maclA _{s}(\cc {2d}),
\label{Eq:PowerKernelSpaces1}
\intertext{and their duals}
&\wideparen \maclA _{0,\flat _\sigma}'(\cc {2d}),
&\quad
&\wideparen \maclA _{\flat _\sigma}'(\cc {2d}),
&\quad
&\wideparen \maclA _{0,s}'(\cc {2d}),
&\quad
&\wideparen \maclA _{s}'(\cc {2d}),
\label{Eq:PowerKernelSpaces2}
\end{alignat}
where
$$
\wideparen \maclA _{\flat _\sigma}(\cc {2d})
=
\sets {K}{(z,w)\mapsto K(z,\overline w)\in \maclA _{\flat _\sigma}(\cc {2d})},
$$
and similarly for the other spaces in \eqref{Eq:PowerKernelSpaces1}
and \eqref{Eq:PowerKernelSpaces2}. In the end we prove that if
$s>0$ or $s=\flat _\sigma$, then the integral operators in \eqref{Eq:CompKernelOp},
\begin{alignat}{4}
T_K\,  &: &\, \maclA _s(\cc d) &\mapsto & \maclA _s'(\cc d)
\quad &\text{when}& \quad
K &\in \wideparen \maclA _s'(\cc {2d}),
\label{Eq:SpecMapIntro1}
\intertext{and}
T_K\,  &: &\, \maclA _{0,s}(\cc d) &\mapsto & \maclA _{0,s}'(\cc d)
\quad &\text{when}& \quad
K &\in \wideparen \maclA _{0,s}'(\cc {2d}),
\label{Eq:SpecMapIntro2}
\end{alignat}
are uniquely defined and continuous, and similarly when the roles
of the non-duals in \eqref{Eq:AnalSpaces1} and
\eqref{Eq:PowerKernelSpaces1}, and their duals in \eqref{Eq:AnalSpaces2}
and \eqref{Eq:PowerKernelSpaces2} are swapped.
We also prove the opposite direction, that any linear
and continuous operators between such spaces are
given by such kernel operators. These kernel results
are given in Propositions \ref{Prop:KernelsDifficultDirection} and
\ref{Prop:KernelsEasyDirection}. Due to
the Bargmann transform homeomorphisms, these results
are also equivalent to Theorems 3.3 and 3.4 in \cite{CheSigTo2}
on kernel theorems for Pilipovi{\'c} spaces. (See Subsection \ref{subsec1.5}.)

\par

Note that, if $s\ge \frac 12$, then the spaces of power series expansions
above can be identified with certain spaces of analytic and semi-conjugate analytic
functions. For example we have
\begin{alignat*}{2}
\maclA _{\flat _\sigma}(\cc {d})
&=
\sets {F\in A(\cc {d})}{|F(z)|\lesssim
e^{r|z|^{\frac {2\sigma}{\sigma +1}}}\ \text{for some}\ r>0},
\quad \sigma >0,
\\[1ex]
\maclA _s(\cc {d})
&=
\sets {F\in A(\cc {d})}{|F(z)|\lesssim
e^{\frac 12\cdot |z|^2-r|z|^{\frac 1{2s}}}\ \text{for some}\ r>0}
\\[1ex]
\maclA _s'(\cc {d})
&=
\sets {F\in A(\cc {d})}{|F(z)|\lesssim
e^{\frac 12\cdot |z|^2+r|z|^{\frac 1{2s}}}\ \text{for every}\ r>0}
\\[1ex]
\maclA _{\flat _\sigma}'(\cc {d})
&=
\sets {F\in A(\cc {d})}{|F(z)|\lesssim
e^{r|z|^{\frac {2\sigma}{\sigma -1}}}\ \text{for every}\ r>0},
\quad \sigma >1,
\\[1ex]
\maclA _{\flat _1}'(\cc {d}) &= A(\cc {d})
\quad \text{and}\quad
\maclA _{0,\flat _1}'(\cc {d}) = A_d(\{ 0\} ),
\end{alignat*}
and similarly for $\wideparen \maclA _s(\cc {2d})$
and $\wideparen \maclA _s'(\cc {2d})$. In particular, the mappings
\eqref{Eq:SpecMapIntro1} and \eqref{Eq:SpecMapIntro2} can be formulated
in terms of those function spaces.

\par

If instead $s\in (0,\frac 12]$ and $\sigma >0$,
and $t\in \mathbf C$, then $K(z,w)\mapsto K(z,w)e^{t(z,w)}$ is homeomorphic on
$$
\wideparen A_{\flat _\sigma}'(\cc {2d}),
\quad
\wideparen A_{0,\flat _\sigma}'(\cc {2d}),
\quad
\wideparen A_{s}'(\cc {2d})
\quad \text{and}\quad
\wideparen A_{0,s}'(\cc {2d}),
$$
see Theorem \ref{Thm:TtMap}.
In particular,
\eqref{Eq:CompPseudo2} implies that
the mappings \eqref{Eq:SpecMapIntro1} and \eqref{Eq:SpecMapIntro2}
still hold true with $\op _{\mathfrak V}(a)$ in place of $T_K$. (Cf. Theorems \ref{Thm:KernelsPseudoDifficultDirection} and \ref{Thm:KernelsPseudoEasyDirection}.)

\par

In the case $s\ge \frac 12$,
the conditions on $a$ and its kernel
$K_a$ of $\op _{\mathfrak V}(a)$ are slightly different. More precisely,
these conditions are of the form
\begin{align*}
|a(z,w)| &\lesssim e^{\frac 12\cdot |z-w|^2+r(|z|^{\frac 1{2s}} + |w|^{\frac 1{2s}})}
\intertext{and}
|K(z,w)| &\lesssim e^{\frac 12\cdot (|z|^2+|w|^2)+r(|z|^{\frac 1{2s}} + |w|^{\frac 1{2s}})}
\end{align*}
in order for the mappings \eqref{Eq:SpecMapIntro1} and \eqref{Eq:SpecMapIntro2}
should hold. (Cf. Theorems \ref{Thm:GSPseudoDifficultDirection} and
\ref{Thm:GSPseudoEasyDirection}.)

\par

In Section \ref{sec3} we consider operators \eqref{Eq:CompKernelOp},
where certain linear pullbacks of their kernels obey suitable mixed and weighted
Lebesgue norm estimates. We prove that such operators are continuous between
appropriate  (weighted) Lebesgue spaces of entire functions. For example, let $\omega$
be a weight on $\cc d\times \cc d$ and $\omega _1,\omega _2$ be weights
on $\cc d$ such that
$$
\frac {\omega _2(z)}{\omega _1(w)}\lesssim \omega (z,w)
$$
and let
$$
G_{K,\omega}(z,w) = K_\omega (z,z+w),
$$
where
$$
K_\omega (z,w)=e^{-\frac 12(|z|^2+|w|^2)}|K(z,w)|\omega (\sqrt 2 \overline z,\sqrt 2 w).
$$
If $p,q,p_j,q_j\in [1,\infty ]$ satisfy
$$
\frac 1{p_1}-\frac 1{p_2}=\frac 1{q_1}-\frac 1{q_2}=1-\frac 1{p}-\frac 1{q}
\quad \text{and}\quad
q\le p,
$$
and $G_{K,\omega}\in L^{p,q}(\cc d\times \cc d)$, then it follows from
Theorem \ref{Thm:LebOpCont} that $T_K$ is continuous
from  $A^{\mabfp _1}_{E,(\omega _1)}(\cc d)$ to $A^{\mabfp _2}_{E,(\omega _2)}(\cc d)$. 
By slightly modifying the definition of $G_{K,\omega}$ we also deduce another similar but
different continuity result where the condition $q\le p$ above is removed (cf.
Theorem \ref{Thm:LebOpCont3}).

\par

We also present some consequences of these results.
Theorem \ref{Thm:LebOpCont2} can be considered as a special case of
Theorem \ref{Thm:LebOpCont} formulated by analytic pseudo-differential
operators instead of integral operators. Theorems
\ref{Thm:PseudoModCont} and \ref{Thm:PseudoModCont2} are obtained
by imposing conditions on moderateness on $\omega$, $\omega _1$ and $\omega _2$
above and translating Theorem \ref{Thm:LebOpCont} and  \ref{Thm:LebOpCont3}
to real pseudo-differential operators via the Bargmann transform and its inverse.
These approaches show that obtained continuity results on analytic pseudo-differential
or integral operators might be suitable when investigating real pseudo-differential
operators. In fact, Theorems \ref{Thm:PseudoModCont} and
\ref{Thm:PseudoModCont2} agree with 
the sharp results
\cite[Theorem 3.3]{To10}, \cite[Theorem 3.1]{Toft19} and \cite[Theorem 2.2]{Toft22}
in the Banach space case. Remark \ref{Rem:Some Extensions} in the end of Section
\ref{sec3} shows that
our approach can be used to extend the latter results on real pseudo-differential operators
to include situations
with non-moderate weights. We note that the moderate
condition on weights may in some situations be significantly restrictive (cf.
Remark \ref{Rem:BasicModSpaceProperties} in Section \ref{sec1}).

\par

\section{Preliminaries}\label{sec1}

\par

In this section we recall some facts on involved function and distribution
spaces as well as on pseudo-differential operators. In Subsection
\ref{subsec1.1} we introduce suitable weight classes. Thereafter we recall in
Subsections \ref{subsec1.2}--\ref{subsec1.4} the definitions and basic
properties for Gelfand-Shilov, Pilipovi{\'c} and modulation spaces. Then we
discuss in Subsection \ref{subsec1.5} the Bargmann transform and
recall some topological spaces of entire
functions or power series expansions on $\cc d$. The section is concluded with a
review of some facts on pseudo-differential operators.

\par

\subsection{Weight functions}\label{subsec1.1}
A \emph{weight} on $\rr d$ is a positive function $\omega \in  L^\infty _{loc}(\rr d)$
such that $1/\omega \in  L^\infty _{loc}(\rr d)$. The weight $\omega$ on $\rr d$
is called \emph{moderate} if there is a positive locally bounded function
$v$ on $\rr d$ such that
\begin{equation}\label{eq:2}
\omega(x+y)\le C\omega(x)v(y),\quad x,y\in\rr{d},
\end{equation}
for some constant $C\ge 1$. If $\omega$ and $v$ are weights on $\rr d$ such
that \eqref{eq:2} holds, then $\omega$ is also called \emph{$v$-moderate}.
The set of all moderate weights on $\rr d$ is denoted by $\mascP _E(\rr d)$.

\par

The weight $v$ on $\rr d$ is called \emph{submultiplicative},
if it is even and \eqref{eq:2}
holds for $\omega =v$. From now on, $v$ always denotes
a submultiplicative
weight if nothing else is stated. In particular,
if \eqref{eq:2} holds and $v$ is submultiplicative, then it follows
by straight-forward computations that
\begin{equation}\label{eq:2Next}
\begin{gathered}
\frac {\omega (x)}{v(y)} \lesssim \omega(x+y) \lesssim \omega(x)v(y),
\\[1ex]
\quad
v(x+y) \lesssim v(x)v(y)
\quad \text{and}\quad v(x)=v(-x),
\quad x,y\in\rr{d}.
\end{gathered}
\end{equation}
Here and in what follows we write
$A(\theta )\lesssim B(\theta )$, $\theta \in \Omega$,
if there is a constant $c>0$ such that $A(\theta )\le cB(\theta )$
for all $\theta \in \Omega$.

\par

If $\omega$ is a moderate weight on $\rr d$, then by \cite{To11}
and above, there is a
submultiplicative weight
$v$ on $\rr d$ such that \eqref{eq:2} and \eqref{eq:2Next}
hold (see also \cite{Groch,To11}). Moreover if $v$ is
submultiplicative on $\rr d$, then
\begin{equation}\label{Eq:CondSubWeights}
1\lesssim v(x) \lesssim e^{r|x|}
\end{equation}
for some constant $r>0$ (cf. \cite{Groch}). In particular, if $\omega$ is moderate, then
\begin{equation}\label{Eq:ModWeightProp}
\omega (x+y)\lesssim \omega (x)e^{r|y|}
\quad \text{and}\quad
e^{-r|x|}\le \omega (x)\lesssim e^{r|x|},\quad
x,y\in \rr d
\end{equation}
for some $r>0$.

\par

\subsection{Gelfand-Shilov spaces}\label{subsec1.2}
Let $0<s \in \mathbf R$ be fixed. Then the (Fourier invariant)
Gelfand-Shilov space $\maclS _s(\rr d)$ ($\Sigma _s(\rr d)$) of
Roumieu type (Beurling type) consists of all $f\in C^\infty (\rr d)$
such that
\begin{equation}\label{gfseminorm}
\nm f{\maclS _{s,h}^\sigma}\equiv \sup \frac {|x^\alpha \partial ^\beta
f(x)|}{h^{|\alpha  + \beta |}(\alpha !\, \beta !)^s}
\end{equation}
is finite for some $h>0$ (for every $h>0$). Here the supremum should be taken
over all $\alpha ,\beta \in \mathbf N^d$ and $x\in \rr d$. The semi-norms
$\nm \cdo {\maclS _{s,h}^\sigma}$ induce an inductive limit topology for the
space $\maclS _s(\rr d)$ and projective limit topology for $\Sigma _s(\rr d)$, and
the latter space becomes a Fr{\'e}chet space under this topology.

\par

The space $\maclS _s(\rr d)\neq \{ 0\}$ ($\Sigma _s(\rr d)\neq \{0\}$), if and only if
$s\ge \frac 12$ ($s> \frac 12$).

\par

The \emph{Gelfand-Shilov distribution spaces} $\maclS _s'(\rr d)$
and $\Sigma _s'(\rr d)$ are the dual spaces of $\maclS _s(\rr d)$
and $\Sigma _s(\rr d)$, respectively.

\par

We have
\begin{equation}\label{GSembeddings}
\begin{aligned}
\maclS _{1/2} (\rr d) &\hookrightarrow \Sigma _s  (\rr d) \hookrightarrow
\maclS _s (\rr d)
\hookrightarrow  \Sigma _t(\rr d)
\\[1ex]
&\hookrightarrow
\mascS (\rr d)
\hookrightarrow \mascS '(\rr d) 
\hookrightarrow \Sigma _t' (\rr d)
\\[1ex]
&\hookrightarrow  \maclS _s'(\rr d)
\hookrightarrow  \Sigma _s'(\rr d) \hookrightarrow \maclS _{1/2} '(\rr d),
\quad \frac 12<s<t.
\end{aligned}
\end{equation}
Here and
in what follows we use the notation $A\hookrightarrow B$ when the topological
spaces $A$ and $B$ satisfy $A\subseteq B$ with continuous embeddings.

\par

A convenient family of functions concerns the Hermite functions
$$
h_\alpha (x) = \pi ^{-\frac d4}(-1)^{|\alpha |}
(2^{|\alpha |}\alpha !)^{-\frac 12}e^{\frac {|x|^2}2}
(\partial ^\alpha e^{-|x|^2}),\quad \alpha \in \nn d.
$$
The set of Hermite functions on $\rr d$ is an orthonormal basis for
$L^2(\rr d)$. It is also a basis for the Schwartz space and its distribution space,
and for any $\Sigma _s$ when $s>\frac 12$,
$\maclS _s$ when $s\ge \frac 12$ and their distribution
spaces. They are also eigenfunctions to the Harmonic
oscillator $H=H_d\equiv |x|^2-\Delta$ and to the Fourier transform
$\mathscr F$, given by
$$
(\mathscr Ff)(\xi )= \widehat f(\xi ) \equiv (2\pi )^{-\frac d2}\int _{\rr
{d}} f(x)e^{-i\scal  x\xi }\, dx, \quad \xi \in \rr d,
$$
when $f\in L^1(\rr d)$. Here $\scal \cdo \cdo$ denotes the usual
scalar product on $\rr d$. In fact, we have
$$
H_dh_\alpha = (2|\alpha |+d)h_\alpha .
$$

\par

The Fourier transform $\mathscr F$ extends
uniquely to homeomorphisms on $\mathscr S'(\rr d)$,
$\maclS _s'(\rr d)$ and on $\Sigma _s'(\rr d)$. Furthermore,
$\mascF$ restricts to
homeomorphisms on $\mathscr S(\rr d)$,
$\maclS _s(\rr d)$ and on $\Sigma _s (\rr d)$,
and to a unitary operator on $L^2(\rr d)$. Similar facts hold true
when the Fourier transform is replaced by a partial
Fourier transform.

\par

Gelfand-Shilov spaces and their distribution spaces can also
be characterized by estimates of short-time Fourier
transform, (see e.{\,}g. \cite{GZ,Teof,Toft18}).
More precisely, let $\phi \in \mascS  (\rr d)$ be
fixed.
Then the \emph{short-time
Fourier transform} $V_\phi f$ of $f\in \mascS '
(\rr d)$ with respect to the \emph{window function} $\phi$ is
the Gelfand-Shilov distribution on $\rr {2d}$, defined by
$$
V_\phi f(x,\xi )  =
\mascF (f \, \overline {\phi (\cdo -x)})(\xi ), \quad x,\xi \in \rr d.
$$
If $f ,\phi \in \mascS (\rr d)$, then it follows that
$$
V_\phi f(x,\xi ) = (2\pi )^{-\frac d2}\int _{\rr d} f(y)\overline {\phi
(y-x)}e^{-i\scal y\xi}\, dy, \quad x,\xi \in \rr d.
$$

\par

By \cite[Theorem 2.3]{To11} it follows that the definition of the map
$(f,\phi)\mapsto V_{\phi} f$ from $\mascS (\rr d) \times \mascS (\rr d)$
to $\mascS(\rr {2d})$ is uniquely extendable to a continuous map from
$\maclS _s'(\rr d)\times \maclS_s'(\rr d)$
to $\maclS_s'(\rr {2d})$, and restricts to a continuous map
from $\maclS _s (\rr d)\times \maclS _s (\rr d)$
to $\maclS _s(\rr {2d})$.
The same conclusion holds with $\Sigma _s$ in place of
$\maclS_s$, at each place.

\par

In the following propositions we give characterizations of Gelfand-Shilov
spaces and their distribution spaces in terms of estimates of the short-time Fourier transform.
We omit the proof since the first part follows from \cite[Theorem 2.7]{GZ})
and the second part from \cite[Proposition 2.2]{Toft18}.
See also \cite{CPRT10} for related results.

\par

\begin{prop}\label{stftGelfand2}
Let $s\ge \frac 12$ ($s>\frac 12$), $\phi \in \maclS _s(\rr d)\setminus 0$
($\phi \in \Sigma _s(\rr d)\setminus 0$) and let $f$ be a
Gelfand-Shilov distribution on $\rr d$. Then the following is true:
\begin{enumerate}
\item $f\in \maclS _s (\rr d)$ ($f\in \Sigma_s(\rr d)$), if and only if
\begin{equation}\label{stftexpest2}
|V_\phi f(x,\xi )| \lesssim  e^{-r (|x|^{\frac 1s}+|\xi |^{\frac 1s})}, \quad x,\xi \in \rr d,
\end{equation}
for some $r > 0$ (for every $r>0$).
\item $f\in \maclS _s'(\rr d)$ ($f\in \Sigma _s'(\rd)$), if and only if
\begin{equation}\label{stftexpest2Dist}
|V_\phi f(x,\xi )| \lesssim  e^{r(|x|^{\frac 1s}+|\xi |^{\frac 1\sigma})}, \quad
x,\xi \in \rr d,
\end{equation}
for every $r > 0$ (for some $r > 0$).
\end{enumerate}
\end{prop}

\par

\subsection{Spaces of Hermite series and power series expansions}
\label{subsec1.3}

\par

Next we recall the definitions of topological vector spaces of
Hermite series expansions, given in \cite{Toft18}. As in \cite{Toft18},
it is convenient to use suitable extensions of
$\mathbf R_+$ when indexing our spaces.

\par

\begin{defn}
The sets $\mathbf R_\flat$ and $\overline {\mathbf R_\flat}$ are given by
$$
{\textstyle{\mathbf R_\flat = \mathbf R_+ \underset{\sigma >0}{\textstyle{\bigcup}}
\{ \flat _\sigma \} }}
\quad \text{and}\quad
{\textstyle{\overline {\mathbf R_\flat} = \mathbf R_\flat \bigcup \{ 0 \} }}.
$$

\par

Moreover, beside the usual ordering in $\mathbf R$, the elements $\flat _\sigma$
in $\mathbf R_\flat$ and $\overline {\mathbf R_\flat}$ are ordered by
the relations $x_1<\flat _{\sigma _1}<\flat _{\sigma _2}<x_2$, when
$\sigma _1$, $\sigma _2$, $x_1$ and $x_2$ are positive real numbers such that
$x_1<\frac 12$ and $x_2\ge \frac 12$.
\end{defn}

\par

\begin{defn}\label{DefSeqSpaces}
Let $p\in [1,\infty ]$, $s\in {\mathbf R_\flat}$, $r\in \mathbf R$, $\vartheta$
be a weight on $\nn d$, and let
$$
\vartheta _{r,s}(\alpha )\equiv
\begin{cases}
e^{r|\alpha |^{\frac 1{2s}}}, & \text{when}\quad s\in \mathbf R_+,
\\[1ex]
r^{|\alpha |}(\alpha !)^{\frac 1{2\sigma}}, & \text{when}\quad s = \flat _\sigma ,
\quad \qquad \alpha \in \nn d.
\end{cases}
$$
Then,
\begin{enumerate}
\item $\ell _0' (\nn d)$ is the set of all sequences $\{c_\alpha \} _{\alpha \in \nn d}
\subseteq \mathbf C$ on $\nn d$;

\vrum

\item $\ell _{0,0}(\nn d)\equiv \{ 0\}$, and $\ell _0(\nn d)$ is the set of all sequences
$\{c_\alpha \} _{\alpha \in \nn d}\subseteq \mathbf C$ such that $c_\alpha \neq 0$
for at most finite numbers of $\alpha$;

\vrum

\item $\ell ^p_{[\vartheta ]}(\nn d)$ is the Banach space which consists of
all sequences $\{ c_\alpha \} _{\alpha \in \nn d} \subseteq \mathbf C$
such that
$$
\nm {\{ c_\alpha \} _{\alpha \in \nn d} }{\ell ^p_{[\vartheta ]}}\equiv
\nm {\{ c_\alpha \vartheta (\alpha )\} _{\alpha \in \nn d} }{\ell ^p} < \infty ;
$$

\vrum

\item $\ell _{0,s}(\nn d)\equiv \underset {r>0}\bigcap \ell ^p_{[\vartheta _{r,s}]}(\nn d)$
and $\ell _s(\nn d)\equiv \underset {r>0}\bigcup \ell ^p_{[\vartheta _{r,s}]}(\nn d)$, with
projective respective inductive limit topologies of $\ell ^p_{[\vartheta _{r,s}]}(\nn d)$
with respect to $r>0$;

\vrum

\item $\ell _{0,s}'(\nn d)\equiv \underset {r>0}\bigcup
\ell ^p_{[1/\vartheta _{r,s}]}(\nn d)$ and $\ell _s'(\nn d)\equiv \underset {r>0}
\bigcap \ell ^p_{[1/\vartheta _{r,s}]}(\nn d)$, with
inductive respective projective limit topologies of $\ell ^p_{[1/\vartheta _{r,s}]}(\nn d)$
with respect to $r>0$.
\end{enumerate}
\end{defn}

\par

Let $p\in [1,\infty ]$, and let $\Omega _N$ be the set of all
$\alpha \in \nn d$ such that $|\alpha |\le N$. Then the
topology of $\ell _0(\nn d)$ is defined by the inductive
limit topology of the sets
$$
\Sets {\{ c_\alpha \} _{\alpha \in \nn d} \in \ell _0'(\nn d)}{c_\alpha =0\
\text{when}\ \alpha \notin \Omega _N}
$$
with respect to $N\ge 0$, and whose topology
is given through the semi-norms
\begin{equation}\label{SemiNormsEllSpaces}
\{ c_\alpha \} _{\alpha \in \nn d}\mapsto \nm {\{ c_\alpha \}
_{|\alpha |\le N} }{\ell ^p(\Omega _N)},
\end{equation}
It is clear that these topologies are independent of $p$.
Furthermore, the topology of
$\ell _0' (\nn d)$ is defined by the semi-norms \eqref{SemiNormsEllSpaces}.
It follows that $\ell _0'(\nn d)$ is a Fr{\'e}chet space, and that the topology
is independent of $p$.

\par

Next we introduce spaces of formal Hermite series expansions
\begin{alignat}{2}
f&=\sum _{\alpha \in \nn d}c_\alpha h_\alpha ,&\quad \{ c_\alpha \}
_{\alpha \in \nn d} &\in \ell _0' (\nn d).\label{Hermiteseries}
\intertext{and power series expansions}
F&=\sum _{\alpha \in \nn d}c_\alpha e_\alpha ,&\quad \{c_\alpha \}
_{\alpha \in \nn d} &\in \ell _0' (\nn d).\label{Powerseries}
\end{alignat}
which correspond to
\begin{equation}\label{ellSpaces}
\ell _{0,s}(\nn d),\quad \ell _s(\nn d),
\quad \ell _s'(\nn d)\quad \text{and}\quad \ell _{0,s}'(\nn d).
\end{equation}
Here
\begin{equation} \label{Eq:basiselements}
e_\alpha (z) \equiv \frac {z^\alpha}{\sqrt {\alpha !}},\qquad z\in \cc d,\ \alpha \in \nn d.
\end{equation}
We consider the mappings
\begin{equation}
\label{T12Map}
T_{\maclH} : \,
\{ c_\alpha \} _{\alpha \in \nn d} \mapsto \sum _{\alpha \in \nn d}
c_\alpha h_\alpha
\quad \text{and}\quad
T_{\maclA} : \,
\{ c_\alpha \} _{\alpha \in \nn d} \mapsto \sum _{\alpha \in \nn d}
c_\alpha e_\alpha
\end{equation}
between sequences, and formal Hermite series and power series expansions.

\par

\begin{defn}\label{DefclHclASpaces}
If $s\in \overline{\mathbf R_\flat}$, then
\begin{alignat}{2}
\maclH _{0,s}(\rr d),\quad \maclH _s(\rr d),
\quad \maclH _s'(\rr d)\quad \text{and}\quad \maclH _{0,s}'(\rr d),\label{clHSpaces}
\intertext{and}
\maclA _{0,s}(\cc d),\quad \maclA _s(\cc d),
\quad \maclA _s'(\cc d)\quad \text{and}\quad \maclA _{0,s}'(\cc d),\label{clASpaces}
\end{alignat}
are the images of $T_{\maclH}$ and $T_{\maclA}$ respectively in \eqref{T12Map}
of corresponding spaces in \eqref{ellSpaces}.
The topologies of the spaces in \eqref{clHSpaces} and \eqref{clASpaces}
are inherited from the corresponding spaces in \eqref{ellSpaces}.
\end{defn}

\par

Since locally absolutely convergent power series expansions can be identified with
entire functions, several of the spaces in \eqref{clASpaces} are identified
with topological vector spaces contained in $A(\cc d)$ (see Theorem
\ref{Thm:AnalSpacesChar} below and the introduction). Here $A(\Omega _0)$
is the set of all (complex valued) functions which are analytic in $\Omega _0$.
(For $\Omega _0\subseteq \cc d$, $A(\Omega _0)=\bigcup A(\Omega )$,
where the union is taken over all open $\Omega \subseteq \cc d$ which
contain $\Omega _0$.)

\par

We recall that $f\in \mascS (\rr d)$ if and only if it can be written as
\eqref{Hermiteseries} such that
$$
|c_\alpha |\lesssim \eabs \alpha ^{-N},
$$
for every $N\ge 0$ (cf. e.{\,}g. \cite{RS}).
In particular it follows from the definitions that the
inclusions
\begin{multline}\label{inclHermExpSpaces}
\maclH _0(\rr d)\hookrightarrow \maclH _{0,s}(\rr d)\hookrightarrow
\maclH _{s}(\rr d) \hookrightarrow \maclH _{0,t}(\rr d)
\\[1ex]
\hookrightarrow \mascS (\rr d) \hookrightarrow \mascS '(\rr d)
\hookrightarrow
\maclH _{0,t}'(\rr d) \hookrightarrow \maclH _{s}'(\rr d)
\\[1ex]
\hookrightarrow \maclH _{0,s}'(\rr d)\hookrightarrow \maclH _0'(\rr d),
\quad \text{when}\ s,t\in \mathbf R_\flat ,\ s<t,
\end{multline}
are dense.

\par

\begin{rem}\label{Rem:LinkEllHs}
By the definition it follows that $T_{\maclH}$ in \eqref{T12Map} is a homeomorphism
between any of the spaces in \eqref{ellSpaces} and corresponding space
in \eqref{clHSpaces}, and that $T_{\maclA}$ in \eqref{T12Map} is a homeomorphism
between any of the spaces in \eqref{ellSpaces} and corresponding space
in \eqref{clASpaces}.
\end{rem}

\par

The next results give some characterizations of $\maclH _s(\rr d)$ and
$\maclH _{0,s}(\rr d)$ when $s$ is a non-negative real number.

\par

\begin{prop}\label{Prop:PilSpacChar1}
Let $0\le s\in \mathbf R$ and let $f\in \maclH _0'(\rr d)$.
Then $f\in \maclH _s(\rr d)$ ($f\in \maclH _{0,s}(\rr d)$), if and only if
$f\in C^\infty (\rr d)$ and satisfies
\begin{equation}\label{GFHarmCond}
\nm{H_d^Nf}{L^\infty}\lesssim h^NN!^{2s},
\end{equation}
for some $h>0$ (every $h>0$). Moreover, it holds
\begin{align*}
\maclH _s(\rr d) &= \maclS _s(\rr d) \neq \{ 0\},
\quad
\maclH _{0,s}(\rr d) = \Sigma _s(\rr d) \neq \{ 0\} 
\quad \text{when}
\quad s\in (\textstyle{\frac 12},\infty),
\\[1ex]
\maclH _s(\rr d) &= \maclS _s(\rr d) \neq \{ 0\},
\quad
\maclH _{0,s}(\rr d) \neq \Sigma _s(\rr d) = \{ 0\} 
\quad \text{when}
\quad s=\textstyle{\frac 12},
\\[1ex]
\maclH _s(\rr d) &\neq \maclS _s(\rr d) = \{ 0\},
\quad
\maclH _{0,s}(\rr d) \neq \Sigma _s(\rr d) = \{ 0\} 
\quad \text{when}
\quad s\in (0,\textstyle{\frac 12}),
\\[1ex]
\maclH _s(\rr d) &\neq \maclS _s(\rr d) = \{ 0\},
\quad
\maclH _{0,s}(\rr d) = \Sigma _s(\rr d) = \{ 0\} 
\quad \text{when}
\quad s=0.
\end{align*}
\end{prop}

\par

We refer to \cite{Toft18} for the proof of Proposition \ref{Prop:PilSpacChar1}.

\par

Due to the pioneering investigations related
to Proposition \ref{Prop:PilSpacChar1} by Pilipovi{\'c} in
\cite{Pil1,Pil2}, we call the spaces $\maclH _s(\rr d)$ and $\maclH _{0,s}(\rr d)$
\emph{Pilipovi{\'c} spaces of Roumieu and Beurling types}, respectively.
In fact, in the restricted case $s\ge \frac 12$,
Proposition \ref{Prop:PilSpacChar1} was proved
already in \cite{Pil1,Pil2}.

\medspace

Later on it will also be convenient for us to have the following definition. Here we let
$F(z_2,\overline z_1)$ and $F(\overline z_2,z_1)$ be the formal power series
\begin{equation}\label{Eq:FormalConjPowerSeries}
\sum
c(\alpha _2,\alpha _1)e_{\alpha _2}(z_2)e_{\alpha _1}(\overline z_1)
\quad \text{and}\quad
\sum
c(\alpha _2,\alpha _1)e_{\alpha _2}(\overline z_2)e_{\alpha _1}(z_1),
\end{equation}
respectively, when $F(z_2,z_1)$ is the formal power series
\begin{equation}\label{Eq:FormalNonConjPowerSeries}
\sum
c(\alpha _2,\alpha _1)e_{\alpha _2}(z_2)e_{\alpha _1}(z_1).
\end{equation}
Here $z_j \in \cc {d_j}$, $j=1,2$, and the sums should be taken over all
$(\alpha _2,\alpha _1)\in \nn {d_2}
\times \nn {d_1}$.

\par

\begin{defn}\label{Def:tauSpaces}
Let $d=d_2+d_1$, $s\in \overline{\mathbf R_\flat}$, $\Theta _{C,1}$
and $\Theta _{C,2}$ be the operators
$$
(\Theta _{C,1}F)(z_2,z_1)=F(z_2,\overline z_1)
\quad \text{and}\quad
(\Theta _{C,2}F)(z_2,z_1)=F(\overline z_2,z_1) 
$$
between formal power series in \eqref{Eq:FormalConjPowerSeries}
and \eqref{Eq:FormalNonConjPowerSeries}, $z_j \in \cc {d_j}$, $j=1,2$.
Then
\begin{equation}\label{Eq:SesAnalSp}
\wideparen \maclA _{0,s}(\cc {d_2}\times \cc {d_1}),
\quad
\wideparen \maclA _s(\cc {d_2}\times \cc {d_1}),
\quad
\wideparen \maclA _s'(\cc {d_2}\times \cc {d_1}),
\quad
\wideparen \maclA _{0,s}'(\cc {d_2}\times \cc {d_1})
\end{equation}
are the images of \eqref{clASpaces} under $\Theta _{C,1}$, and
$\wideparen A (\cc {d_2}\times \cc {d_1})$ and
$\wideparen A _{d_2,d_1}(\{ 0\})=\wideparen A _{d_2+d_1}(\{ 0\})$
are the images of $A(\cc d)$ and
$A _{d_2+d_1}(\{ 0\} )$ respectively under $\Theta _{C,1}$.
The topologies of the spaces in \eqref{Eq:SesAnalSp},
$\wideparen A (\cc {d_2}\times \cc {d_1})$ and
$\wideparen A _{d_2,d_1}(\{ 0\} )$ are inherited from the topologies
in the spaces \eqref{clASpaces}, $A(\cc d)$ and
$A _d(\{ 0\} )$, respectively.
\end{defn}

\par

\begin{rem}\label{Rem:SpaceSpecCase}
By letting $d_2=d$ and $d_1=0$, it follows that $A(\cc d)$
and the spaces in \eqref{clASpaces} can be considered as special cases of
$\wideparen A (\cc {d_2}\times \cc {d_1})$ and the spaces
in \eqref{Eq:SesAnalSp}.

\par

Since $\maclA _{\flat _1}'(\cc d) = A(\cc d)$ and $\maclA _{0,\flat _1}'(\cc d)
= A_d(\{ 0\} )$, it follows that
\begin{equation}\label{Eq:AnalIdents}
\begin{aligned}
\wideparen
\maclA _{\flat _1}'(\cc {d_2}\times \cc {d_1})
&=
\wideparen A(\cc {d_2}\times \cc {d_1}),
\\[1ex]
\wideparen
\maclA _{0,\flat _1}'(\cc {d_2}\times \cc {d_1})
&=
\wideparen A _{d_2,d_1}(\{ 0 \} ).
\end{aligned}
\end{equation}
\end{rem}

\par

The following results are now immediate consequences of Theorems
4.1, 4.2, 5.2 and 5.3 in \cite{Toft18} and Definition \ref{Def:tauSpaces}.
Here let
\begin{align}
\kappa _{1,r,s}(z)
&=
\begin{cases}
e^{r(\log \eabs z)^{\frac 1{1-2s}}}, & s<\frac 12
\\[1ex]
e^{r|z|^{\frac {2\sigma}{\sigma +1}}}, & s=\flat _\sigma ,\ \sigma >0,
\\[1ex]
e^{\frac {|z|^2}2-r|z|^{\frac 1s}}, & s\ge \frac 12,
\end{cases}
\label{Eq:kappa1Def}
\intertext{and}
\kappa _{2,r,s}(z)
&=
\begin{cases}
e^{r|z|^{\frac {2\sigma}{\sigma -1}}}, & s=\flat _\sigma ,\ \sigma >1,
\\[1ex]
e^{\frac {|z|^2}2+r|z|^{\frac 1s}}, & s\ge \frac 12,
\end{cases}
\label{Eq:kappa2Def}
\end{align}

\par

\begin{thm}\label{Thm:AnalSpacesChar}
Let $s_1,s_2\in \mathbf R_\flat$ be such that $s_2>\flat _1$, and let $\kappa _{1,r,s}$
and $\kappa _{2,r,s}$ be given by
\eqref{Eq:kappa1Def} and \eqref{Eq:kappa2Def} respectively, when $r>0$.
Then the following is true:
\begin{enumerate}
\item $\wideparen \maclA _{s_1} (\cc {d_2}\times \cc {d_1})$
($\wideparen \maclA _{0,s_1} (\cc {d_2}\times \cc {d_1})$) consists
of all $K\in \wideparen A(\cc {d_2}\times \cc {d_1})$ such that
$|K|\lesssim \kappa _{1,r,s_1}$ for some $r>0$ (for every $0<r<\frac 12$).

\vrum

\item $\wideparen \maclA _{s_2}' (\cc {d_2}\times \cc {d_1})$
($\wideparen \maclA _{0,s_2}'(\cc {d_2}\times \cc {d_1})$) consists
of all $K\in \wideparen A(\cc {d_2}\times \cc {d_1})$ such that
$|K|\lesssim \kappa _{2,r,s_2}$ for every $r>0$ (for some $r>0$).
\end{enumerate}
\end{thm}

\par

By Remark \ref{Rem:SpaceSpecCase} it follows that Theorem
\ref{Thm:AnalSpacesChar} remains true after the spaces in
\eqref{Eq:SesAnalSp} are replaced by corresponding spaces
in \eqref{clASpaces}.

\par

\subsection{Modulation spaces}\label{subsec1.4}

\par

Before giving the definition of a broad family of modulation spaces,
we make a review of mixed normed spaces of Lebesgue types, adapted to
suitable bases of the Euclidean space $\rr d$.
Let $E$ be the ordered basis $\{ e_1,\dots,e_d \}$
of $\rr d$. Then the ordered basis $E' =\{
e'_1,\dots ,e'_d\}$ (the dual basis of $E$) satisfies
$$
\scal {e_j} {e'_k} = 2\pi \delta_{jk}
\quad \text{for every}\quad
j,k =1,\dots, d.
$$
The corresponding parallelepiped, lattice,
dual parallelepiped and dual lattice are given by
\begin{align*}
\kappa (E) &= \sets{x_1e_1+\cdots+x_de_d}{(x_1,\dots,x_d)
\in\rr d,\ 0\leq
x_k\leq 1,\ k=1,\dots,d},
\\[1ex]
\Lambda _E &=\sets{j_1e_1+\cdots +j_de_d}{(j_1,\dots,j_d)\in
\zz d},
\\[1ex]
\kappa (E') &=\sets {\xi _1e'_1+\cdots+\xi _de'_d}{(\xi _1,\dots ,\xi _d)
\in \rr d,
\ 0\leq \xi _k\leq 1,\ k=1,\dots ,d},
\intertext{and}
\Lambda'_E &= \Lambda_{E'}=\sets{\iota _1e'_1+\cdots +\iota _de'_d}
{(\iota _1,\dots ,\iota _d) \in \zz d},
\end{align*}
respectively. Note here that the Fourier
analysis with respect to general biorthogonal bases has recently been
developed in \cite{RuTo}.

\par

We observe that there is a matrix $T_E$ such that
$e_1,\dots ,e_d$ and $e_1',\dots ,e_d'$ are the images of
the standard basis under $T_E$ and
$T_{E'}= 2\pi(T^{-1}_E)^t$, respectively.

\par

In the following we let
$$
\max (\mabfq ) =\max (q_1,\dots ,q_d)
\quad \text{and}\quad
\min (\mabfq ) =\min (q_1,\dots ,q_d)
$$
when $\mabfq =(q_1,\dots ,q_d)\in [1,\infty ]^d$.

\par

\begin{defn}\label{Def:MixedLebSpaces}
Let $E$ be an ordered basis of $\rr d$ and
$\mabfp =(p_1,\dots ,p_d)\in [1,\infty ]^{d}$.
If  $f\in L^1_{loc}(\rr d)$, then $\nm f{L^{\mabfp }_{E}}$ is
defined by
$$
\nm f{L^{\mabfp }_{E}}\equiv
\nm {g_{d-1}}{L^{p_{d}}(\mathbf R)}
$$
where  $g_k(z_k)$, $z_k\in \rr {d-k}$,
$k=0,\dots ,d-1$, are inductively defined as
\begin{align*}
g_0(x_1,\dots ,x_{d}) &\equiv |f(x_1e_1+\cdots +x_{d}e_d)| ,\quad
(x_1,\dots,x_d) \in \rr d, 
\intertext{and}
g_k(z_k) &\equiv
\nm {g_{k-1}(\cdo ,z_k)}{L^{p_k}(\mathbf R)},
\quad z_k\in \rr {d-k},\ k=1,\dots ,d-1.
\end{align*}
The space $L^{\mabfp }_{E}(\rr d)$ consists
of all $f\in L^1_{loc}(\rr d)$ such that
$\nm f{L^{\mabfp}_{E}}$ is finite, and is called
\emph{$E$-split Lebesgue space (with respect to $\mabfp$)}.
\end{defn}

\par

Next we discuss suitable conditions for bases in the phase space
$\rr {2d}$. We let $\upsigma (X,Y)$ be the
standard symplectic form on the phase space, given by
$$
\upsigma (X,Y) = \scal y\xi -\scal x\eta ,
\qquad
X=(x,\xi )\in \rr {2d},\ Y=(y,\eta )\in \rr {2d}.
$$
We notice that if
\begin{equation}\label{Eq:SymplBasis}
\{ e_1,\dots ,e_d,\ep _1,\dots ,\ep _d\}
\end{equation}
is the standard basis of $\rr {2d}$, then
\begin{equation}\label{Eq:SymplBasisRel}
\upsigma (e_j,e_k) = 0,\quad
\upsigma (e_j,\ep _k) = -\delta _{j,k},
\quad \text{and}\quad
\upsigma (\ep _{j},\ep _{k}) = 0,
\end{equation}
when $j,k\in \{1,\dots ,d\}$. More generally, a basis in \eqref{Eq:SymplBasis}
for the phase space $\rr {2d}$ is called \emph{symplectic} if
\eqref{Eq:SymplBasisRel} holds. A symplectic basis \eqref{Eq:SymplBasis}
for $\rr {2d}$ is called \emph{phase split} if $e_1,\dots ,e_d$
and $\ep _1,\dots ,\ep _d$ span
$$
\sets {(x,0)\in \rr {2d}}{x\in \rr d}
\quad \text{and}\quad
\sets {(0,\xi )\in \rr {2d}}{\xi \in \rr d},
$$
respectively.

\par

Next we give the definition of our class of modulation spaces.

\par

\begin{defn}\label{Def:ModSpaces}
Let $E$ be an ordered basis for $\rr {2d}$, $\mabfp \in [1,\infty ]^{2d}$,
$\phi (x)=\pi ^{-\frac d4}e^{-\frac 12 \cdot |x|^2}$ and let $\omega$ be a weight on $\rr {2d}$.
Then the modulation space
$M^{\mabfp}_{E,(\omega )}(\rr d)$ consists of all $f\in \maclH _{\flat _1}'(\rr d)$
such that
\begin{equation}\label{Eq:ModNormGen}
\nm f{M^{\mabfp}_{E,(\omega )}}\equiv \nm {V_\phi f \cdot \omega}{L^{\mabfp}_E}
\end{equation}
is finite.
\end{defn}

\par

We remark that if $\phi (x)=\pi ^{-\frac d4}e^{-\frac 12 \cdot |x|^2}$ and
$f\in \maclH _{\flat _1}'(\rr d)$, then $(x,\xi )\mapsto V_\phi f(x,\xi )$
is a smooth function (cf. \cite{Toft18}). Furthermore, by \cite[Theorem 4.8]{Toft18} we
get the following. The proof is omitted.

\par

\begin{prop}\label{Prop:ModBanach}
Let $E$ be an ordered basis for $\rr {2d}$, $\mabfp \in [1,\infty ]^{2d}$ and
let $\omega$ be a weight on $\rr {2d}$. Then $M^{\mabfp}_{E,(\omega )}(\rr d)$
is a Banach space with norm given by \eqref{Eq:ModNormGen}.
\end{prop}

\par

If the weight  $\omega$ in Definition \ref{Def:ModSpaces} is a moderate weight,
then we can say more concerning $M^{\mabfp}_{E,(\omega )}(\rr d)$.
In what follows
we let $p'\in [1,\infty ]$ be the conjugate exponent of $p\in [1,\infty ]$, i.{\,}e.
$\frac 1p+\frac 1{p'}=1$.

\par

\begin{prop}\label{Prop:BasicModSpaceProperties}
Let $E$ be an ordered basis for $\rr {2d}$, $\mabfp \in [1,\infty ]^{2d}$ and let
$\omega ,v\in \mascP _E(\rr {2d})$ be such that $\omega$ is $v$-moderate. Then
the following is true:
\begin{enumerate}
\item $\Sigma _1(\rr d)\hookrightarrow M^{\mabfp}_{E,(\omega )}(\rr d) \hookrightarrow
\Sigma _1'(\rr d)$. If in addition $\max (\mabfp )<\infty$, then $\Sigma _1(\rr d)$ is dense in
$M^{\mabfp}_{E,(\omega )}(\rr d)$;

\vrum

\item if $\phi \in M^1_{(v)}(\rr d)\setminus \{0\}$ and $f\in \Sigma _1'(\rr d)$, then
$f\in M^{\mabfp}_{E,(\omega )}(\rr d)$, if and only if the right-hand side of
\eqref{Eq:ModNormGen} is finite. Furthermore, different choices of
$\phi \in M^1_{(v)}(\rr d)\setminus \{0\}$ in \eqref{Eq:ModNormGen}
give rise to equivalent norms;

\vrum

\item $M^{\mabfp}_{E,(\omega )}(\rr d)$ increases with $p_1,\dots ,p_{2d}$ and decreases
with $\omega$;

\vrum

\item if $\mabfp ' = (p_1',\dots ,p_{2d}')$, then the restriction of the $L^2(\rr d)$
scalar product $(\cdo ,\cdo )$ to $\Sigma _1(\rr d)$ is uniquely extendable to a (semi-conjugate)
duality between $M^{\mabfp}_{E,(\omega )}(\rr d)$ and
$M^{\mabfp '}_{E,(1/\omega )}(\rr d)$. If in addition $\max (\mabfp )<\infty$, then
the dual of $M^{\mabfp}_{E,(\omega )}(\rr d)$ can be identified by
$M^{\mabfp '}_{E,(1/\omega )}(\rr d)$ through the form $(\cdo ,\cdo )$.
\end{enumerate}
\end{prop}

\par

Proposition \ref{Prop:BasicModSpaceProperties} follows by similar arguments as in
Chapters 11 and 12 in \cite{Gc2} (see also \cite{To11,Toft18}).

\par

\begin{rem}\label{Rem:BasicModSpaceProperties}
In some sense, the variable $x$ at the weight $\omega (x,\xi )$ in the definition
of modulation spaces quantify growth and decay properties for the involved
functions or distributions. In the same way the variable $\xi$ quantify
regularity or lack of regularity for the involved functions or distributions.

\par

By the analysis in \cite{Toft18} it follows that there are no bounds on
how fast $V_\phi f$ may grow or decay at infinity when
$\phi (x)=\pi ^{-\frac d4}e^{-\frac 12\cdot |x|^2}$ is fixed, $x\in \rr d$,
and $f$ is taken in the class $\maclH _{\flat _1}'(\rr d)$. Since
weights in $\mascP _E(\rr {2d})$ are bounded by exponential functions,
the restrictions of the weights in Proposition \ref{Prop:BasicModSpaceProperties}
are significantly stronger compared to what is the case in Proposition
\ref{Prop:ModBanach}.
A question here concerns wether it is possible to extend
parts of Proposition \ref{Prop:BasicModSpaceProperties} to larger
weight classes than $\mascP _E(\rr {2d})$ or not.

\par

It seems that the invariance properties (2) in Proposition
\ref{Prop:BasicModSpaceProperties} concerning the choice of weight function
are not possible for weights that are not moderate. On the other hand,
(1) and (4) in Proposition \ref{Prop:BasicModSpaceProperties}
hold true for certain
weights outside $\mascP _E(\rr {2d})$. In fact, in \cite{To11},
certain weight classes which contain $\mascP _E(\rr {2d})$
as well as weights of the form
$$
\omega (x,\xi )= e^{\pm r(|x|^{\frac 1s}+|\xi |^{\frac 1s})}, \quad x,\xi \in \rr d,
$$
when $r>0$ and $s>\frac 12$ are introduced.
For corresponding (broader)
families of modulation spaces it is then proved
that Proposition \ref{Prop:BasicModSpaceProperties}
(1) and (4) hold true (with some modifications).
\end{rem}

\par

\subsection{Bargmann transform and spaces of analytic
functions}\label{subsec1.5}

\par

The Bargmann transform $\mathfrak V_d$ is the homeomorphism from
the spaces in \eqref{clHSpaces} to respective spaces in \eqref{clASpaces},
given by $T_{\maclA}\circ T_{\maclH}^{-1}$, where $T_{\maclH}$ and $T_{\maclA}$ are given by
\eqref{T12Map}.


\par

We notice that if $f\in L^p(\rr d)$ for some $p\in [1,\infty ]$,
then $\mathfrak V_df$ is the entire function given by
\begin{equation*}
(\mathfrak V_df)(z) =\pi ^{-d/4}\int _{\rr d}\exp \Big ( -\frac 12(\scal
z z+|y|^2)+2^{1/2}\scal zy \Big )f(y)\, dy,\quad z \in \cc d,
\end{equation*}
which can also be formulated as
$$
(\mathfrak V_df)(z) =\int_{\rr d} \mathfrak A_d(z,y)f(y)\, dy,
\quad z \in \cc d,
$$
or
\begin{equation}\label{bargdistrform}
(\mathfrak V_df)(z) =\scal f{\mathfrak A_d(z,\cdo )},
\quad z \in \cc d,
\end{equation}
where the Bargmann kernel $\mathfrak A_d$ is given by
$$
\mathfrak A_d(z,y)=\pi ^{-d/4} \exp \Big ( -\frac 12(\scal
zz+|y|^2)+2^{1/2}\scal zy\Big ), \quad z \in \cc d, y \in \rr d.
$$
(Cf. \cite{B1,B2}.)
Here
$$
\scal zw = \sum _{j=1}^dz_jw_j\quad \text{and} \quad
(z,w)= \scal z{\overline w}
$$
when
$$
z=(z_1,\dots ,z_d) \in \cc d\quad  \text{and} \quad w=(w_1,\dots ,w_d)\in \cc d,
$$
and otherwise $\scal \cdo \cdo $ denotes the duality between test function
spaces and their corresponding duals which is clear form the context.
We note that the right-hand side in \eqref{bargdistrform} makes sense
when $f\in \maclS _{1/2}'(\rr d)$ and defines an element in $A(\cc d)$,
since $y\mapsto \mathfrak A_d(z,y)$ can be interpreted as an element
in $\maclS _{1/2} (\rr d)$ with values in $A(\cc d)$.

\par

It was proved by Bargmann that $f\mapsto \mathfrak V_df$ is a bijective and isometric map
from $L^2(\rr d)$ to the Hilbert space $A^2(\cc d)$, the set of entire functions $F$ on $\cc
d$ which fullfils
\begin{equation}\label{A2norm}
\nm F{A^2}\equiv \Big ( \int _{\cc d}|F(z)|^2d\mu (z)  \Big )^{1/2}<\infty .
\end{equation}
Recall, $d\mu (z)=\pi ^{-d} e^{-|z|^2}\, d\lambda (z)$, where $d\lambda (z)$ is
the Lebesgue measure on $\cc d$, and the scalar product on $A^2(\cc d)$ is given by
\begin{equation}\label{A2scalar}
(F,G)_{A^2}\equiv  \int _{\cc d} F(z)\overline {G(z)}\, d\mu (z),\quad F,G\in A^2(\cc d).
\end{equation}
For future references we note that the latter scalar product induces the bilinear form
\begin{equation}\label{A2scalarBil}
(F,G)\mapsto \scal FG _{A^2}=\scal FG _{A^2(\cc d)}\equiv
\int _{\cc d} F(z)G(z)\, d\mu (z)
\end{equation}
on $A^2(\cc d)\times \overline{A^2(\cc d)}$.

\par

In \cite{B1} it was proved that the orthonormal basis
$\{ h_\alpha \}_{\alpha \in \nn d}$ in $L^2(\rr d)$ of Hermite
functions is mapped to the orthonormal basis
$\{ e_\alpha \} _{\alpha \in \nn d}$ in $A^2(\cc d)$
(cf. \eqref{Eq:basiselements}).
Furthermore, there is a convenient reproducing formula on
$A^2(\cc d)$. In fact, let $\Pi _A$ be the operator from $L^2(d\mu )$
to $A(\cc d)$, given by
\begin{equation}\label{reproducing}
(\Pi _AF)(z)= \int _{\cc d} F(w)e^{(z,w)}\, d\mu (w),\quad z \in \cc d.
\end{equation}
 Then it is proved in \cite{B1} that $\Pi _A$ is an orthonormal
projection from $L^2(d\mu)$ to $A^2(\cc d)$.

\par

From now on we assume that $\phi$ in the definition of the short-time
Fourier transform is given by
\begin{equation}\label{phidef}
\phi (x)=\pi ^{-d/4}e^{-|x|^2/2}, \quad x\in \rr d,
\end{equation}
if nothing else is stated. For such $\phi$, it follows by straight-forward
computations that the relationship between
the Bargmann transform and the short-time Fourier transform
is given by
\begin{equation}\label{bargstft1}
\mathfrak V_d = U_{\mathfrak V}\circ V_\phi ,\quad \text{and}\quad
U_{\mathfrak V}^{-1} \circ \mathfrak V_d =  V_\phi ,
\end{equation}
where $U_{\mathfrak V}$ is the linear, continuous and bijective operator on
$\mathscr D'(\rr {2d})\simeq \mathscr D'(\cc d)$, given by
\begin{equation}\label{UVdef}
(U_{\mathfrak V}F)(x+i\xi ) = (2\pi )^{d/2} e^{(|x|^2+|\xi |^2)/2}e^{-i\scal x\xi}
F(2^{1/2}x,-2^{1/2}\xi ), \quad x,\xi \in \rr d,
\end{equation}
cf. \cite{To11}.

\par

\begin{defn}\label{thespaces}
Let $E$ be an ordered basis for $\rr {2d}$, $U_{\mathfrak V}$ be the
operator in \eqref{UVdef}, $\mabfp \in [1,\infty ]^{2d}$,
$\phi (x)=\pi ^{-\frac d4}e^{-\frac 12 \cdot |x|^2}$ and let $\omega$
be a weight on $\rr {2d}$.
\begin{enumerate}
\item The space $B^{\mabfp}_{E,(\omega )}(\cc d)$ consists of all
$F\in L^1_{loc}(\cc d)$ such that
$$
\nm F{B^{\mabfp}_{E,(\omega )}}\equiv
\nm {(U_{\mathfrak V}^{-1}F)\cdot \omega }{L^{\mabfp}_E}
$$
is finite;

\vrum

\item The space $A^{\mabfp}_{E,(\omega )}(\cc d)$ consists of all $F\in A(\cc
d)\bigcap B^{\mabfp}_{E,(\omega )}(\cc d)$ with topology inherited
from $B^{\mabfp}_{E,(\omega )}(\cc d)$.
\end{enumerate}
\end{defn}

\par

We note that the spaces in Definition \ref{thespaces} are normed
spaces when $\min (\mabfp) \ge 1$.

\par

For conveneincy we set $\nm F{B^{\mabfp}_{E,(\omega )}}=\infty$, when
$F\notin B^{\mabfp}_{E,(\omega )}(\cc d)$ is measurable, and
$\nm F{A^{\mabfp}_{E,(\omega )}}=\infty$, when $F\in A(\cc d)\setminus
B^{\mabfp}_{E,(\omega )}(\cc d)$.

\par

\begin{rem}\label{Rem:SpaceNotations}
In Definitions \ref{Def:ModSpaces} and \ref{thespaces},
important cases appear when $E$ is the standard basis
for $\rr {2d}$ and
$p_1=\cdots =p_d=p\in [1,\infty ]$ and $p_{d+1}
=\cdots =p_{2d}=q\in [1,\infty]$.
For such choices of $E$ and $\mabfp$ we set
$L^{p,q} = L^{\mabfp}_E$,
\begin{alignat*}{4}
M^{p,q}_{(\omega )} &= M^{\mabfp}_{E,(\omega )},&
\quad
A^{p,q}_{(\omega )} &= A^{\mabfp}_{E,(\omega )} &
\quad &\text{and} & \quad
B^{p,q}_{(\omega )} &= B^{\mabfp}_{E,(\omega )}.
\intertext{We also set}
M^p_{(\omega )} &= M^{p,p}_{(\omega )},&
\quad
A^p_{(\omega )} &=A^{p,p}_{(\omega )} &
\quad &\text{and} & \quad
B^p_{(\omega )} &=B^{p,p}_{(\omega )}.
\intertext{If in addition $\omega =1$, then we set}
M^{\mabfp}_{E,(\omega )} &= M^{\mabfp}_{E}, &
\quad
M^{p,q}_{(\omega )} &= M^{p,q} &
\quad &\text{and} &\quad
M^p_{(\omega )} &= M^p,
\intertext{and similarly for $A^{\mabfp}_{E,(\omega )}$ and
$B^{\mabfp}_{E,(\omega )}$ spaces.
\newline
\indent
If instead $E =\{ e_{d+1},\dots ,e_{2d},e_1,\dots ,e_d \}$
where $e_1,\dots ,e_{2d}$ is the standard basis for $\rr {2d}$
and
$p_1=\cdots =p_d=q\in [1,\infty ]$ and $p_{d+1}
=\cdots =p_{2d}=p\in [1,\infty]$, then we set $L^{p,q}_*=L^{\mabfp}_E$,}
W^{p,q}_{(\omega )} &= M^{\mabfp}_{E,(\omega )}, &
\quad
A^{p,q}_{*,(\omega )} &= A^{\mabfp}_{E,(\omega )} &
\quad &\text{and} &\quad
B^{p,q}_{*,(\omega )}&= B^{\mabfp}_{E,(\omega )}.
\end{alignat*}
\end{rem}

\par

We notice that the space $W^{p,q}_{(\omega )}$ in Remark
\ref{Rem:SpaceNotations} is an example of a (weighted)
Wiener amalgam space (cf. \cite{Fe0,F1}).

\par

For future references we observe that the $B^{p}_{(\omega )}$ norm is given by
\begin{multline}\label{BLpnorm}
\nm F{B^{p}_{(\omega )}} = 2^{d/p}(2\pi )^{-d/2}\left ( \int _{\cc d} |e^{-|z|^2/2}
F(z)\omega (2^
{1/2}\overline z)|^p\, d\lambda (z)  \right )^{1/p}
\\[1ex]
= 2^{d/p}(2\pi )^{-d/2}\left ( {{\iint} }_{\! \! \rr {2d}} |e^{-(|x|^2+|\xi |^2)/2}F(x+i\xi )\omega (2^
{1/2}x,-2^{1/2}\xi )|^p\, dxd\xi \right )^{1/p}
\end{multline}
(with obvious modifications when $p=\infty$). Especially it follows that the norm and scalar
product in $B^2_{(\omega )}(\cc d)$ take the forms
\begin{equation*}
\begin{alignedat}{2}
\nm F{B^{2}_{(\omega )}} &= \left ( \int _{\cc {d}} |F(z)\omega (2^
{1/2}\overline z)|^2\, d\mu (z) \right )^{1/2},& \quad F &\in B^2_{(\omega )}(\cc d)
\\[1ex]
(F,G)_{B^2_{(\omega )}} &=  \int _{\cc {d}} F(z)\overline {G(z)}\omega (2^
{1/2}\overline z)^2\, d\mu (z) ,& \quad F,G &\in B^2_{(\omega )}(\cc d)
\end{alignedat}
\end{equation*}
(cf. \eqref{A2norm} and \eqref{A2scalar}).

\par

By the definitions and \eqref{bargstft1}
it follows that the Bargmann transform
is an isometric injection from $M^{\mabfp}_{E,(\omega )}(\rr d)$
to $A^{\mabfp}_{E,(\omega )}(\cc d)$. In fact, we have
the following refinement. We omit the proof since the result
is a special case of Theorem 4.8 in \cite{Toft18}.

\par

\begin{prop}\label{mainstep1prop}
Let $E$ be an ordered basis for $\rr {2d}$, $\mabfp \in [1,\infty ]^{2d}$,
and $\omega$ be a weight on $\rr {2d}$. Then the Bargmann transform
is an isometric bijection from $M^{\mabfp}_{E,(\omega )}(\rr d)$
to $A^{\mabfp}_{E,(\omega )}(\cc d)$.
\end{prop}

\par

%
%
%

Finally, the \emph{SCB transform} (i.{\,}e. the Semi Conjugated
Bargmann transform), $\mathfrak V_{\Theta ,d_2,d_1}$ is
defined as $\Theta _{C,1} \circ \mathfrak V_{d_2+d_1}$. We also set
$\mathfrak V_{\Theta ,d} = \mathfrak V_{\Theta ,d,d}$.  Evidently, all properties
of the Bargmann transform
carry over to analogous properties for the SCB transform. Assume that
$E$ is a basis for $\rr {2d_2}\times \rr {2d_1}$,
$\mabfp \in [1,\infty ]^{2d_2+2d_1}$, $p,q\in (0,\infty ]$ and that
$\omega$ is a weight on $\rr {2d_2}\times \rr {2d_1}$
Then
$\wideparen A^{\mabfp}_{E,(\omega )}(\cc {d_2+d_1})$
is the image of $A^{\mabfp}_{E,(\omega )}(\cc {d_2+d_1})$
under the map $\Theta _{C,1}$ with the topology
defined by the norm
$$
\nm {a}{\wideparen A^{\mabfp}_{E,(\omega )}}\equiv
\nm {\Theta _{C,1} a}{A^{\mabfp}_{E,(\Theta _{C,1} \omega )}},
\qquad a\in \wideparen A^{\mabfp}_{E,(\omega )}(\cc {d_2+d_1}).
$$
The spaces
$$
\wideparen A^{p,q}_{(\omega )}(\cc {d_2}\times \cc {d_1}),
\quad
\wideparen A^{p}_{(\omega )}(\cc {d_2}\times \cc {d_1}),
\quad
\wideparen A^{p,q}(\cc {d_2}\times \cc {d_1})
\quad \text{and}\quad
\wideparen A^{p}(\cc {d_2}\times \cc {d_1})
$$
and their norms, and the scalar product
$(\cdo ,\cdo )_{\wideparen A^{2}}$ are defined analogously.

\par

\subsection{Pseudo-differential operators}\label{subsec1.6}

\par

Next we recall some properties in pseudo-differential calculus.
Let $\GL (d,\Omega)$ be the set of $d\times d$-matrices with
entries in the set $\Omega$, $a\in \Sigma _1
(\rr {2d})$, and let $A\in \GL (d,\mathbf R)$ be fixed. Then the
pseudo-differential operator $\op _A(a)$
is the linear and continuous operator on $\Sigma _1 (\rr d)$, given by
\begin{equation}\label{e0.5}
(\op _A(a)f)(x)
=
(2\pi  ) ^{-d}\iint a(x-A(x-y),\xi )f(y)e^{i\scal {x-y}\xi }\,
dyd\xi, \quad x\in \rr d.
\end{equation}
For general $a\in \Sigma _1'(\rr {2d})$, the
pseudo-differential operator $\op _A(a)$ is defined as the continuous
operator from $\Sigma _1(\rr d)$ to $\Sigma _1'(\rr d)$ with
distribution kernel
\begin{equation}\label{atkernel}
K_{a,A}(x,y)=(2\pi )^{-d/2}(\mascF _2^{-1}a)(x-A(x-y),x-y), \quad 
x,y \in \rr d.
\end{equation}
Here $\mascF _2F$ is the partial Fourier transform of $F(x,y)\in
\Sigma _1'(\rr {2d})$ with respect to the $y$ variable. This
definition makes sense since the mappings
\begin{equation}\label{homeoF2tmap}
\mascF _2\quad \text{and}\quad F(x,y)\mapsto F(x-A(x-y),x-y)
\end{equation}
are homeomorphisms on $\Sigma _1'(\rr {2d})$.
In particular, the map $a\mapsto K_{a,A}$ is a homeomorphism on
$\Sigma _1'(\rr {2d})$.

\par

The standard (Kohn-Nirenberg) representation, $a(x,D)=\op (a)$, and
the Weyl quantization $\op ^w(a)$ of $a$ are obtained by choosing
$A=0$ and $A=\frac 12 I$, respectively, in \eqref{e0.5} and \eqref{atkernel},
where $I$ is the identity matrix.

\par

\begin{rem}\label{BijKernelsOps}
By Fourier's inversion formula, \eqref{atkernel} and the kernel theorem
\cite[Theorem 2.2]{LozPer}, \cite[Theorem 2.5]{Teof2} for operators from
Gelfand-Shilov spaces to their duals,
it follows that the map $a\mapsto \op _A(a)$ is bijective from $\Sigma _1'(\rr {2d})$
to the set of all linear and continuous operators from $\Sigma _1(\rr d)$
to $\Sigma _1'(\rr {2d})$.
\end{rem}

\par

By Remark \ref{BijKernelsOps}, it follows that for every $a_1\in \Sigma _1'(\rr {2d})$
and $A_1,A_2\in \GL (d,\mathbf R)$, there is a unique $a_2\in \Sigma _1'(\rr {2d})$ such that
$\op _{A_1}(a_1) = \op _{A_2} (a_2)$. By Section 18.5 in \cite{Ho1},
the relation between $a_1$ and $a_2$
is given by
\begin{equation}
\label{calculitransform}
\op _{A_1}(a_1) = \op _{A_2}(a_2)
\quad \Leftrightarrow \quad
a_2=e^{i\scal {(A_1-A_2)D_\xi}{D_x}}a_1.
\end{equation}
Here we note that the operator $e^{i\scal {AD_\xi}{D_x}}$ is homeomorphic
on $\Sigma _1(\rr {2d})$ and its dual
(cf. \cite{CaTo,CarWal,Tr}). For modulation spaces we have the following
subresult of Proposition 2.8 in \cite{Toft22}.

\par

\begin{prop}\label{Prop:ExpOpSTFT}
Let $s\ge \frac 12$, $A\in \GL (d,\mathbf R)$, $p,q\in (0,\infty ]$,
$\phi ,a\in \Sigma _1(\rr {2d})$ and let $T_A = e^{i\scal{AD_\xi}{D_x}}$.
If $\omega \in \mascP _E(\rr {4d})$ and
$$
\omega _A(x,\xi ,\eta ,y) = \omega (x+Ay,\xi +A^*\eta ,\eta ,y),
$$
then $T_A$ from $\Sigma _1(\rr {2d})$ to $\Sigma _1(\rr {2d})$
extends uniquely to a homeomorphism from $M^{p,q}_{(\omega )}(\rr {2d})$
to $M^{p,q}_{(\omega _A)}(\rr {2d})$, and
\begin{equation}\label{Eq:ExpOpModSp}
\nm {T_Aa}{M^{p,q}_{(\omega _A)}} \asymp
\nm a{M^{p,q}_{(\omega )}}.
\end{equation}
\end{prop}

\section{Kernel theorems and analytic pseudo-differential
operators}\label{sec2}

\par

In the first part of the section we show that there is a one to one correspondence
between linear and continuous mappings from $\maclA _s'$ to $\maclA _s$
($\maclA _s$ to $\maclA _s'$) and mappings with kernels in $\wideparen A_s$
($\wideparen A_s'$) with respect to the measure $d\mu$ (cf. Propositions
\ref{Prop:KernelsDifficultDirection} and \ref{Prop:KernelsEasyDirection}).
Thereafter we deduce in Theorems
\ref{Thm:KernelsPseudoDifficultDirection}--\ref{Thm:GSPseudoEasyDirection}
analogous results for analytic pseudo-differential
operators based on Theorem \ref{Thm:TtMap} which deals with mapping
properties of the operator which takes $a(z,w)$ into $e^{(z,w)}a(z,w)$.

\par

Here and in what follows, any extension of the $A^2$-form, $(\cdo ,\cdo )_{A^2}$
from  $\maclA _0(\cc d)\times \maclA _0(\cc d)$ 
to $\mathbf C$ is still called $A^2$-form and still denoted by $(\cdo ,\cdo )_{A^2}$.
Similar approaches yield extensions of the forms $\scal \cdo \cdo _{A^2}$
and $(\cdo ,\cdo )_{\wideparen A^2}$.

\par

By the definitions, $\ell _s'(\nn d)$ and $\ell _{0,s}'(\nn d)$ are the duals of
$\ell _s(\nn d)$ and $\ell _{0,s}(\nn d)$, respectively, through unique extensions
of the $\ell ^2(\nn d)$ form on $\ell _0(\nn d)$. Since the spaces in
\eqref{clASpaces} are images of the spaces in \eqref{ellSpaces}
under the map $T_{\maclA}$ in \eqref{T12Map}, the following lemma is an
immediate consequence of these duality properties. The result is also implicitly
given in \cite{CheSigTo2,Toft18}.

\par

\begin{lemma}\label{Lemma:AnDual}
Let $s\in \overline{\mathbf R_\flat}$. Then the following is true:
\begin{enumerate}
\item the form $(F,G)\mapsto (F,G)_{A^2}$
from  $\maclA _0(\cc d)\times \maclA _0(\cc d)$ to $\mathbf C$ is uniquely extendable
to continuous forms from $\maclA _s(\cc d)\times \maclA _s'(\cc d)$ to $\mathbf C$,
and from $\maclA _{0,s}(\cc d)\times \maclA _{0,s}'(\cc d)$ to $\mathbf C$. Furthermore,
the duals of $\maclA _{s}(\cc d)$ and $\maclA _{0,s}(\cc d)$ can be identified by
$\maclA _{s}'(\cc d)$ and $\maclA _{0,s}'(\cc d)$ through the form $(\cdo ,\cdo )_{A_2}$;

\vrum

\item the form $(F,G)\mapsto \scal FG _{A^2}$
from  $\maclA _0(\cc d)\times \overline{\maclA _0(\cc d)}$ to $\mathbf C$ is
uniquely extendable
to continuous forms from $\maclA _s(\cc d)\times \overline {\maclA _s'(\cc d)}$
to $\mathbf C$,
and from $\maclA _{0,s}(\cc d)\times \overline{\maclA _{0,s}'(\cc d)}$ to
$\mathbf C$. Furthermore,
the duals of $\maclA _{s}(\cc d)$ and $\maclA _{0,s}(\cc d)$ can be identified by
$\maclA _{s}'(\cc d)$ and $\maclA _{0,s}'(\cc d)$ through the form $\scal \cdo \cdo _{A_2}$.
\end{enumerate}
\end{lemma}

\par

The following two propositions follow by applying $\mathfrak V_{\Theta ,d_2,d_1}$ on
Theorem 3.3 and 3.4 in \cite{CheSigTo2}, and using
Lemma \ref{Lemma:AnDual}. The details are left for the reader.

\par

\begin{prop}\label{Prop:KernelsDifficultDirection}
Let $s\in \overline{\mathbf R_\flat}$, and let $T$ be a linear and continuous map
from $\maclA _0(\cc {d_1})$ to $\maclA _0'(\cc {d_2})$.
Then the following is true:
\begin{enumerate}
\item if $T$ is a linear and continuous map from $\maclA _s'(\cc {d_1})$ to
$\maclA _s(\cc {d_2})$, then there is a unique
$K\in \wideparen \maclA _s(\cc{d_2}\times \cc {d_1})$ such that
\begin{equation}\label{Kmap2}
TF =\big ( z_2\mapsto \scal {K(z_2,\cdo)}F _{A^2(\cc {d_1})} \big )
\end{equation}
holds true;

\vrum

\item if $T$ is a linear and continuous map from $\maclA _s(\cc {d_1})$ to
$\maclA _s'(\cc {d_2})$, then there is a unique
$K\in \wideparen \maclA _s'(\cc {d_2}\times \cc {d_1})$ such that \eqref{Kmap2} holds true.
\end{enumerate}

\par

The same holds true if $\maclA _s$, $\wideparen \maclA _s$, $\maclA _s'$
and $\wideparen \maclA _s'$ are replaced by
$\maclA _{0,s}$, $\wideparen \maclA _{0,s}$, $\maclA _{0,s}'$ and $\wideparen
\maclA _{0,s}'$, respectively, at each occurrence.
\end{prop}

\par

\begin{prop}\label{Prop:KernelsEasyDirection}
Let $K\in \wideparen \maclA _0'(\cc {d_2}\times \cc {d_1})$,
$s\in \overline{\mathbf R_\flat}$ and let $T$ be the linear and
continuous map from $\maclA _0(\cc {d_1})$ to
$\maclA _0'(\cc {d_2})$, given by
\begin{equation}\label{Kmap}
F\mapsto TF = \big ( z_2\mapsto \scal {K(z_2,\cdo)}F _{A^2(\cc {d_1})} \big ).
\end{equation}
Then the following is true:
\begin{enumerate}
\item if $K\in \wideparen \maclA _s(\cc {d_2}\times \cc {d_1})$, then $T$ extends uniquely
to a linear and continuous map from $\maclA _s'(\cc {d_1})$ to
$\maclA _s(\cc {d_2})$;

\vrum

\item if $K\in \wideparen \maclA _s'(\cc {d_2}\times \cc {d_1})$, then $T$ extends uniquely
to a linear and continuous map from $\maclA _s(\cc {d_1})$ to
$\maclA _s'(\cc {d_2})$.
\end{enumerate}

\par

The same holds true if $\maclA _s$, $\wideparen \maclA _s$, $\maclA _s'$
and $\wideparen \maclA _s'$ are replaced by
$\maclA _{0,s}$, $\wideparen \maclA _{0,s}$, $\maclA _{0,s}'$ and $\wideparen
\maclA _{0,s}'$, respectively, at each occurrence.
\end{prop}

\par

The operator $T$ in \eqref{Kmap} should be interpreted as $T$ in
the formula
$$
(TF,G)_{A^2(\cc {d_2})} = (K,G\otimes \overline{F})_{A^2(\cc {d_2}\times \cc {d_1})},
\qquad F\in \maclA _0(\cc {d_1}),\ G\in \maclA _0(\cc {d_2}).
$$

\par

Next we recall the definition of analytic pseudo-differential operators.
(See \cite[Definition 6.20]{To11} in the case $t=0$, as well as
\cite{Bau,Berezin71}.)

\par

\begin{defn}\label{APsiDO}
Let $ a \in \wideparen \maclA _{\flat_1} ' (\cc {d} \times \cc {d} )$. Then the
\emph{analytic pseudo-differential operator}
$\op _{\mathfrak V}(a)$ (A$\Psi$DO) with symbol $a$ is given by
\begin{multline}\label{APDO}
(\op _{\mathfrak V}(a)F)(z)
=
\int_{\cc {d} } a(z,w)F(w)
e^{(z,w)}\, d\mu (w)
\\[1ex]
=
(F, \overline {a(z,\cdo )}e^{(\cdo ,z)})_{A^2(\cc {d})}, \quad
z \in \cc d.
\end{multline}
\end{defn}

\par

By the definition it follows that the relation between the operator
kernel $K$ and the symbol $a$ is given by
$$
K (z,w) = e^{(z,w)} a(z,w),\quad z,w\in \cc {d},
$$
provided the multiplication on the right-hand side makes sense.

\par

This leads to the question about mapping properties of $T_t$ defined by
\begin{equation} \label{T_t}
(T_ta) (z,w) = e^{t(z,w)} a(z,w), \quad
z,w\in \cc {d},\ t \in \mathbf C
\end{equation}
when $a$ belongs to a suitable subspace of
$\wideparen \maclA _s '(\cc {2d})$.

\par

First we notice that $ a \in \wideparen A (\cc {2d})$, if and only if $T_t a
\in \wideparen A (\cc {2d})$, and that the inverse of
$ T_t $ is $T_{-t} $. Hence $ T_t $ is well-defined and a homeomorphism
on $ \wideparen \maclA _{\flat_1} '(\cc {2d}).$

\par

If $T_{\maclH}$ is the same as in \eqref{T12Map} then we shall investigate
the map $T_{0,t} $ in the commutative diagram:
\begin{equation}  \label{commdiagram}
\begin{CD}
\ell _{\flat_1} ' (\nn {2d})    @>T_{0,t}>> \ell _{\flat_1} ' (\nn {2d})
\\
@V T_{\maclH}VV        @VVT_{\maclH}V\\
\wideparen \maclA _{\flat_1}  (\cc {2d})     @>>T_t>
\wideparen \maclA _{\flat_1}  (\cc {2d}).
\end{CD}
\end{equation}

\par

Therefore, let $ a \in \wideparen \maclA _{\flat_1} '(\cc {2d})$
with the expansion
$$
a(z,w) = \sum_{\alpha, \beta \in \nn {d}}
c(\alpha ,\beta ) e_{\alpha }(z)e_{\beta}(\overline w) =
\sum_{\alpha, \beta \in \nn {d}} c(\alpha ,\beta )
\frac{z^{\alpha} \overline w ^\beta}{\sqrt{\alpha! \beta !}},
\quad z,w\in \cc d,
$$
where
\begin{equation} \label{coeffestimate}
|c(\alpha ,\beta )| \lesssim r^{|\alpha + \beta|}\sqrt{\alpha! \beta !}
\end{equation}
for every $r>0.$
Since
$$
e^{t(z,w)} = \sum_{\gamma \in \nn {d}}
\frac{t^{|\gamma|}z ^{\gamma} \overline w ^\gamma}{\gamma !},
\quad z,w\in \cc {d},
$$
we have
\begin{equation} \label{e-times-symbol}
e^{t(z,w)}  a(z,w) = \sum_{\gamma \in \nn {d}} \sum_{\alpha, \beta \in \nn {d}}
\varphi_{t,z,w} (\alpha, \beta, \gamma),
\quad z,w\in \cc {d},\ t \in \mathbf C,
\end{equation}
where
\begin{equation} \label{varphi}
\varphi_{t,z,w} (\alpha, \beta, \gamma) =
\frac{c(\alpha ,\beta) t^{|\gamma|}z ^{\alpha + \gamma} \overline w
^{\beta + \gamma}}{\gamma !\sqrt{\alpha! \beta !}},
\quad z,w\in \cc {d},\ t \in \mathbf C,
\end{equation}

\par

We shall prove that the series in \eqref{e-times-symbol} is locally
uniformly convergent with respect to $t,z$ and $w$.
If $ |t| < R$, $ |z| < R $ and $ |w| < R $ for some fixed $R>0$, then
by \eqref{coeffestimate} we get
$$
|\varphi_{t,z,w} (\alpha, \beta, \gamma) | \lesssim
\frac{r^{|\alpha + \beta|} R^{|\alpha + \beta + 2\gamma|}}{\gamma !}
\le
\frac{(rR)^{|\alpha + \beta|} (dR^2)^{|\gamma|}}{|\gamma| !}
$$
for all $ \alpha, \beta, \gamma \in \nn {d}$.
Since the series
$$
\sum_{\alpha, \beta, \gamma \in \nn {d}}\frac{(rR)^{|\alpha + \beta|}
(2R^2)^{|\gamma|}}{|\gamma| !}
$$
is convergent when $r$ is chosen strictly smaller than $ R^{-1}$,
the asserted uniform convergence follows from Weierstass' theorem.

\par

In particular, we may change the order of summation in
\eqref{e-times-symbol} to obtain
\begin{multline*}
(T_t a) (z,w) =  \sum_{\alpha, \beta \in \nn {d}}\sum_{\gamma \in \nn {d}}
\frac{c(\alpha ,\beta) t^{|\gamma|}z ^{\alpha + \gamma}
\overline w ^{\beta + \gamma}}{\gamma !\sqrt{\alpha! \beta !}}
\\[1ex]
=  \sum_{\alpha, \beta \in \nn {d}}\sum_{\gamma \in \nn {d}}
c(\alpha ,\beta) t^{|\gamma|}
\left ( { {\alpha + \gamma} \choose \gamma }
{ {\beta + \gamma} \choose \gamma } \right )^{1/2}
e_{\alpha + \gamma }(z) e_{\beta + \gamma}(\overline w)
\\[1ex]
=  \sum_{\alpha, \beta \in \nn {d}} (T_{0,t} c) (\alpha ,\beta) e_{\alpha }(z)
e_{\beta }(\overline w),
\quad z,w\in \cc {d},\ t \in \mathbf C,
\end{multline*}
where
\begin{equation} \label{T_{0,t} c}
(T_{0,t} c) (\alpha ,\beta) = \sum_{\gamma \leq \alpha ,\beta}
c(\alpha - \gamma, \beta - \gamma) t^{|\gamma|}
\left ( { \alpha  \choose \gamma } { \beta  \choose \gamma } \right )^{1/2},
\quad t \in \mathbf C,
\end{equation}
and we have identified $T_{0,t} $ in the diagram \eqref{commdiagram}.

\par

We have now the following:

\par

\begin{prop}\label{Prop:T_{0,t} homeomorphism}
Let $t \in \cc {},$ $ s,s_0 \in \overline{\mathbf R_\flat}$ be such that
$ s < 1/2 $ and $ 0 < s_0 \leq 1/2$, and let
$T_{0,t} $ be the map on $ \ell _0' (\nn {2d}) $  given by \eqref{T_{0,t} c}.
Then $T_{0,t} $ is a continuous and bijective map  on $ \ell _0' (\nn {2d}) $
with the inverse  $T_{0,-t} $. Furthermore,  $T_{0,t} $ restricts to
homeomorphism from $ \ell _s ' (\nn {2d}) $ to $ \ell _s ' (\nn {2d}) $, and
from $ \ell _{0,s_0} ' (\nn {2d}) $ to $ \ell _{0,s_0} ' (\nn {2d}) $.
\end{prop}

\par

\begin{proof}
The topology on $ \ell _0' (\nn {2d}) $ can be defined by the family of semi-norms
$$
p_N (\{ c(\alpha, \beta)\}_{\alpha,\beta \in \nn {d}})
= \sup_{|\alpha| \leq N}\sup_{|\beta| \leq N} |c(\alpha, \beta )|, \quad
N \in \nn {}.
$$
Then, for a given $c\in  \ell _0' (\nn {2d}) $ we have
\begin{multline*}
p_N (T_{0,t}( c)) \leq \sup_{|\alpha|, |\beta| \leq N}
 \sum_{\gamma \leq \alpha ,\beta} c(\alpha - \gamma, \beta - \gamma)
 |t|^{|\gamma|}
\left ( { \alpha  \choose \gamma } { \beta  \choose \gamma } \right )^{1/2}
\\[1ex]
\leq p_N (c) \sum_{|\gamma| \leq N}  |t|^{|\gamma|} 2^N
\\[1ex]
\leq 2^N (1+|t|)^{Nd} p_N (c), \quad t\in \cc {}, N \in \nn {},
\end{multline*}
and the continuity of $T_{0,t} $  on $ \ell _0' (\nn {2d}) $  follows.

\par

By straight-forward computations it also follows that $T_{0,-t}$ is
the inverse of $T_{0,t}$, which gives asserted homeomorphism
properties of $T_{0,t}$ on $ \ell _0' (\nn {2d}) $.

\par

Next we consider the case when $s,s_0 \in \mathbf R_{+}$. Assume that
$$
|c(\alpha, \beta)| \lesssim e^{\frac 1h({|\alpha|^{1/2s}  + |\beta|^{1/2s}}) },
\quad \alpha,\beta \in \nn {d},
$$
for some constant $h > 0.$
Then
\begin{equation*}
(T_{0,t} c) (\alpha, \beta) \lesssim  \sum_{\gamma \leq \alpha, \beta}
 e^{\frac 1h ({|\alpha-\gamma|^{1/2s}  + |\beta-\gamma|^{1/2s}}) } |t|^{|\gamma|}
\left ( { \alpha  \choose \gamma } { \beta  \choose \gamma } \right )^{1/2}
\le I(\alpha )I(\beta ),
\end{equation*}
where
$$
I(\alpha )= \left ( \sum_{\gamma \leq \alpha} { \alpha  \choose \gamma }
e^{\frac 2h ({ |\alpha-\gamma|^{1/2s}}) } |t|^{|\gamma|}  \right )^{1/2},
\quad t\in \cc {}, \alpha \in \nn {d},
$$
and similarly for $ I(\beta)$. Since
$$
I (\alpha) \leq  e^{\frac 1h{ |\alpha|^{1/2s}} } \left ( \sum_{\gamma \leq \alpha}
{ \alpha  \choose \gamma } |t|^{|\gamma|}  \right )^{1/2}
= e^{\frac 1h{ |\alpha|^{1/2s}} } (1 + |t|)^{d|\alpha|/2}, \quad t\in \cc {}, \alpha \in \nn {d},
$$
we get
$$
|(T_{0,t} c) (\alpha, \beta) |\lesssim
 e^{\frac 1h({ |\alpha|^{1/2s} +  |\beta|^{1/2s}}) } (1 + |t|)^{d(|\alpha|+|\beta|)/2}
 \lesssim
e^{\frac 2h{ (|\alpha|^{1/2s} +  |\beta|^{1/2s})} } ,
$$
$t\in \cc {}, \alpha,\beta \in \nn {d}, $ 
where the last inequality follows from the fact that $ s< 1/2$.
This gives the continuity assertions for $T_{0,t}$ in the case when
$s,s_0 \in \mathbf R_{+}$ and $s,s_0 < 1/2 $.

\par

For $ s_0 = 1/2 $ we have
$$
|(T_{0,t} c) (\alpha, \beta) |\lesssim
e^{\frac1h ( |\alpha| +  |\beta|)} (1 + |t|)^{d(|\alpha|+|\beta|)/2} =
e^{\frac1{h_1}(|\alpha| +  |\beta|) }, \quad
t\in \cc {}, \alpha,\beta \in \nn {d}, 
$$
for some other choice of $h_1 > 0$ which only depend on $|t|,$ $d$
and $h$ and the continuity of  $T_{0,t}$ on $ \ell _{0,\frac{1}{2}} '
(\nn {2d})$ follows.

\par

It remains to consider the case when $ s = s_0 = \flat _{\sigma}$ for
some $ \sigma >0$. Assume that
$$
|c(\alpha, \beta)| \leq C r^{|\alpha + \beta|}
\left  (\alpha! \beta! \right )^{\frac{1}{2\sigma}},\quad
\alpha,\beta \in \nn {d},
$$
for some constants $C,r > 0$. Then
\begin{multline*}
|(T_{0,t} c) (\alpha, \beta)| \leq C  \sum_{\gamma \leq \alpha, \beta}
r^{|\alpha +  \beta - 2\gamma|}\left  ((\alpha-\gamma)!
(\beta-\gamma)! \right )^{\frac{1}{2\sigma}}
|t|^{|\gamma|} \left ( { \alpha  \choose \gamma } { \beta  \choose \gamma }
\right )^{1/2}
\\[1ex]
\leq C_r \sum_{\gamma \leq \alpha, \beta}
r^{|\alpha +  \beta |}\left  (\alpha! \beta! \right )^{\frac{1}{2\sigma}}
|t|^{|\gamma|} 2^{|\alpha +  \beta |/2}
\\[1ex]
\leq C_r (2r)^{|\alpha +  \beta |} \left  (\alpha! \beta! \right )^{\frac{1}{2\sigma}}
\sum_{|\gamma| \leq |\alpha + \beta|} |t|^{|\gamma|}
\\[1ex]
\leq C_r (2r(1+|t|)^d )^{|\alpha +  \beta |} \left  (\alpha! \beta! \right )^{\frac{1}{2\sigma}},
\quad
t\in \cc {}, \alpha,\beta \in \nn {d}, 
\end{multline*}
where $C_r >0$ only depends on $C$ and $r$. This shows that  $T_{0,t}$ is
continuous on $ \ell _{\flat_{\sigma}} ' (\nn {2d})$ and on $ \ell _{{0,\flat_{\sigma}}}
' (\nn {2d})$.
\end{proof}

\par

We have now the following:

\par

\begin{thm}\label{Thm:TtMap}
Let $t \in \mathbf C$, $ s,s_0 \in \overline{\mathbf R_\flat}$ be such that
$ s < \frac 12 $ and $ 0 < s_0 \leq \frac 12$, and let
$T_{t} $ be given by \eqref{T_t}  when $ a \in \wideparen
\maclA _{\flat_1} '(\cc {2d}).$
Then the following is true:
\begin{enumerate}
\item $T_t$ restricts to a homeomorphism from
$ \wideparen \maclA _{0,1/2} '(\cc {2d})$ to
$ \wideparen \maclA _{0,1/2} '(\cc {2d})$;

\vrum

\item  $T_t$ from $ \wideparen \maclA _{0,1/2} '(\cc {2d})$ to
$ \wideparen \maclA _{0,1/2} '(\cc {2d})$ extends uniquely to
homeomorphisms from $ \wideparen \maclA _{s} '(\cc {2d})$
to $ \wideparen \maclA _{s} '(\cc {2d})$ and from
$ \wideparen \maclA _{0,s_0} '(\cc {2d})$ to $ \wideparen \maclA _{0,s_0} '(\cc {2d})$.
\end{enumerate}
\end{thm}

\par

\begin{proof}
By the commutative diagram \eqref{commdiagram} we have
$$
T_t a = (T_1 \circ T_{0,t} \circ T_1 ^{-1}) a,\qquad a \in \wideparen
\maclA _{\flat_1} (\cc {2d}),
$$
and letting
$ T_t a = (T_1 \circ T_{0,t} \circ T_1 ^{-1}) a $
for general $ a \in \wideparen \maclA _{0} '(\cc {2d})$, the continuity
assertions follow from Proposition \ref{Prop:T_{0,t} homeomorphism}.

\par

It remains to prove the uniqueness.
Let $a \in \wideparen \maclA _{0} '(\cc {2d})$,  $b \in \wideparen \maclA _{0} (\cc {2d})$,
with the corresponding expansion coefficients $ c_a (\alpha, \beta)$,
and  $ c_b (\alpha, \beta),$ respectively, and let
$c_{T_t a} (\alpha, \beta)$  be the coefficients of $ T_t a \in
\wideparen \maclA _{0} '(\cc {2d})$, $ \alpha, \beta \in \nn {d}$.
Then
$$
(a,b)_{\wideparen A^2} = \sum_{|\alpha + \beta| \leq N}  c_a (\alpha, \beta)
\overline{ c_b (\alpha, \beta)},
$$
for some $N \in \nn {}$ depending on $b.$ Now choose a sequence
$ a_j \in  \wideparen \maclA _{0} (\cc {2d})$ such that
\begin{equation}\label{Eq:a_jWeakConv}
\lim_{j \rightarrow \infty } (a_j,b)_{\wideparen A^2} = (a,b)_{\wideparen A^2},
\quad \text{for every} \quad b \in \wideparen \maclA _{0} (\cc {2d}).
\end{equation}
If $c_{a_j} (\alpha, \beta)$ denote the coefficients in the expansion of
$a_j, $ $j \in \nn {}$, then it follows from \eqref{Eq:a_jWeakConv}
\begin{equation} \label{Eq:c_jConvergence}
\lim_{j \rightarrow \infty } c_{a_j} (\alpha, \beta) = c_{a} (\alpha, \beta),
\quad \text{for every} \quad \alpha, \beta \in \nn d
\end{equation}
by taking $b(z,w)= e_\alpha (z) e_\beta (\overline w)$.
The uniqueness follows if we prove that
\begin{equation} \label{T_j convergence}
\lim_{j \rightarrow \infty } ((T_t a_j),b)_{\wideparen A^2}
= ((T_t a),b)_{\wideparen A^2},
\quad \text{for every} \quad b \in \wideparen \maclA _{0} (\cc {2d}).
\end{equation}

\par

Let the coefficients of $T_t a_j$ be denoted by $c_{t, a_j}(\alpha, \beta) $,
$\alpha, \beta \in \nn {d}$.
By \eqref{T_{0,t} c} and  \eqref{Eq:c_jConvergence} we get
$c_{t,a_j}(\alpha, \beta) \rightarrow c_{t, a}(\alpha, \beta) $ as
$j\rightarrow \infty$, for every $(\alpha, \beta) \in \nn {2d}$,
and \eqref{T_j convergence} follows since
$$
((T_t a_j),b)_{\wideparen A^2} = \sum_{|\alpha + \beta| \leq N}
c_{t, a_j} (\alpha, \beta)\overline{ c_b (\alpha, \beta)},
$$
and
$$
((T_t a),b)_{\wideparen A^2} = \sum_{|\alpha + \beta| \leq N}
c_{t, a} (\alpha, \beta)\overline{ c_b (\alpha, \beta)},  \quad
b \in \wideparen \maclA _{0} (\cc {2d}),
$$
where $N$ depends on $b$ only.
\end{proof}

\par

The following two theorems now follows by combining Propositions
\ref{Prop:KernelsDifficultDirection} and \ref{Prop:KernelsEasyDirection}
with Theorem \ref{Thm:TtMap}. The details are left for the reader.

\par

\begin{thm}\label{Thm:KernelsPseudoDifficultDirection}
Let $s\in \overline{\mathbf R_\flat}$ be such that $s<\frac 12$ and
let $T$ be a linear and continuous map
from $\maclA _s(\cc {d})$ to $\maclA _s'(\cc {d})$. Then there is a unique
$a\in \wideparen \maclA _s'(\cc{d}\times \cc {d})$ such that
$T =\op _{\mathfrak V}(a)$.

\par

The same holds true if $s<\frac 12$, $\maclA _s$, $\maclA _s'$
and $\wideparen \maclA _s'$ are replaced by $s\le \frac 12$,
$\maclA _{0,s}$, $\maclA _{0,s}'$ and $\wideparen \maclA _{0,s}'$,
respectively, at each occurrence.
\end{thm}

\par

\begin{thm}\label{Thm:KernelsPseudoEasyDirection}
Let $a\in \wideparen \maclA _0'(\cc {d}\times \cc {d})$ and
$s\in \overline{\mathbf R_\flat}$ be such that $s<\frac 12$.
If $a\in \wideparen \maclA _s'(\cc {d}\times \cc {d})$, then $\op _{\mathfrak V}(a)$
extends uniquely to a linear and continuous map from $\maclA _s(\cc {d})$ to
$\maclA _s'(\cc {d})$.

\par

The same holds true if $s<\frac 12$, $\maclA _s$, $\maclA _s'$
and $\wideparen \maclA _s'$ are replaced by $s\le \frac 12$,
$\maclA _{0,s}$, $\maclA _{0,s}'$ and $\wideparen \maclA _{0,s}'$,
respectively, at each occurrence.
\end{thm}

\par

The analogous results to Theorems \ref{Thm:KernelsPseudoDifficultDirection} and \ref{Thm:KernelsPseudoEasyDirection}
for larger $s$  are equivalent to kernel theorems for Fourier invariant
Gelfand-Shilov spaces.

\par

\begin{thm}\label{Thm:GSPseudoDifficultDirection}
Let $s\ge \frac 12$ ($s>\frac 12$). Then the following is true:
\begin{enumerate}
\item If $T$ is a linear and continuous map
from $\maclA _s'(\cc {d})$ to
$\maclA _s(\cc {d})$ (from $\maclA _{0,s}'(\cc {d})$ to
$\maclA _{0,s}(\cc {d})$), then there is a unique
$a\in \wideparen A (\cc{d}\times \cc {d})$ such that
\begin{equation}\label{Eq:GSPseudoSymbCond1}
|a(z,w)|\lesssim e^{\frac 12\cdot |z-w|^2-r(|z|^{\frac 1s}+|w|^{\frac 1s})},
\quad z,w \in \cc {d},
\end{equation}
for some (for every) $r>0$ and $T =\op _{\mathfrak V}(a)$;

\vrum

\item If $T$ is a linear and continuous map
from $\maclA _s(\cc {d})$ to
$\maclA _s'(\cc {d})$ (from $\maclA _{0,s}(\cc {d})$ to
$\maclA _{0,s}'(\cc {d})$), then there is a unique
$a\in \wideparen A (\cc{d}\times \cc {d})$ such that
\begin{equation}\label{Eq:GSPseudoSymbCond2}
|a(z,w)|\lesssim e^{\frac 12\cdot |z-w|^2+r(|z|^{\frac 1s}+|w|^{\frac 1s})},
\quad z,w \in \cc {d},
\end{equation}
for every (for some) $r>0$ and $T =\op _{\mathfrak V}(a)$.
\end{enumerate}
\end{thm}

\par

\begin{thm}\label{Thm:GSPseudoEasyDirection}
Let $s\ge \frac 12$ ($s> \frac 12$). Then the following is true:
\begin{enumerate}
\item If $a\in \wideparen A(\cc {d}\times \cc {d})$ satisfies
\eqref{Eq:GSPseudoSymbCond1} for some (for every) $r>0$,
then $\op _{\mathfrak V}(a)$ from $\maclA _0(\cc d)$ to $\maclA _0'(\cc d)$
is uniquely extendable to a linear and continuous map from
$\maclA _s'(\cc {d})$ to $\maclA _s(\cc {d})$ (from
$\maclA _{0,s}'(\cc {d})$ to $\maclA _{0,s}(\cc {d})$);

\vrum

\item If $a\in \wideparen A(\cc {d}\times \cc {d})$ satisfies
\eqref{Eq:GSPseudoSymbCond2} for every (for some) $r>0$,
then $\op _{\mathfrak V}(a)$ from $\maclA _0(\cc d)$ to $\maclA _0'(\cc d)$
is uniquely extendable to a linear and continuous map from
$\maclA _s(\cc {d})$ to $\maclA _s'(\cc {d})$ (from
$\maclA _{0,s}(\cc {d})$ to $\maclA _{0,s}'(\cc {d})$).
\end{enumerate}
\end{thm}

\par

\begin{proof}
We only prove the results in the Roumieu case. The Beurling case follows by
similar arguments and is left for the reader.

\par

If $T$ is the same as in \eqref{Kmap} for some $K\in \wideparen A(\cc d\times \cc d)$,
then $T=\op _{\mathfrak V}(a)$ when $a(z,w)=e^{-(z,w)}K(z,w)$,
$ (z,w) \in \cc {d} \times \cc {d}$. Since
$$
|e^{-(z,w)}|e^{\frac 12 (|z|^2+|w|^2)} = e^{-\frac 12\cdot |z-w|^2},
\quad z,w \in \cc {d},
$$
Theorem \ref{Thm:AnalSpacesChar} gives
\begin{alignat*}{2}
K\in \wideparen \maclA _s(\cc d\times \cc d)\quad &\Leftrightarrow &
\quad
|K(z,w)| &\lesssim e^{\frac 12\cdot (|z|^2+|w|^2)-r(|z|^{\frac 1s}+|w|^{\frac 1s})}
\\[1ex]
&\Leftrightarrow &
\quad |a(z,w)|
&\lesssim
e^{\frac 12\cdot (|z-w|^2)-r(|z|^{\frac 1s}+|w|^{\frac 1s})},\quad z,w \in \cc {d}.
\end{alignat*}
for some $r>0$. In the same way,
$$
K\in \wideparen \maclA _s'(\cc d\times \cc d)\quad \Leftrightarrow
\quad
|a(z,w)|\lesssim e^{\frac 12\cdot (|z-w|^2)+r(|z|^{\frac 1s}+|w|^{\frac 1s})},
\quad z,w \in \cc {d},
$$
for every $r>0$. The results now follows from these relations and
Propositions \ref{Prop:KernelsDifficultDirection}
and \ref{Prop:KernelsEasyDirection}
\end{proof}

\par

\begin{rem}
For strict subspaces of $\wideparen \maclA _{0,1/2}'(\cc d)$ in Definition
\ref{Def:tauSpaces}, the
estimates imposed on their elements are given by \eqref{Eq:kappa1Def} or
by \eqref{Eq:kappa2Def} for suitable assumptions on $r>0$.
It is evident that in all such
cases, these conditions are violated under the action of $T_t$ in Theorem \ref{Thm:TtMap}
when $t\neq 0$. Hence, Theorem \ref{Thm:TtMap}
cannot be extended to other spaces in Definition \ref{Def:tauSpaces}.

\par

In particular, the conditions \eqref{Eq:GSPseudoSymbCond1} and
\eqref{Eq:GSPseudoSymbCond2} in Theorems
\ref{Thm:GSPseudoDifficultDirection} and
\ref{Thm:GSPseudoEasyDirection} can not be replaced by
the convenient condition that $a$ should belong to e.{\,}g.
\begin{equation}\label{Eq:SpacesInThms}
\wideparen \maclA _{0,s}(\cc {d}\times \cc {d}),
\quad
\wideparen \maclA _s(\cc {d}\times \cc {d}),
\quad
\wideparen \maclA _{s_1}'(\cc {d}\times \cc {d})
\quad \text{or}\quad
\wideparen \maclA _{0,s_2}'(\cc {d}\times \cc {d})
\end{equation}
when $s\in \overline {\mathbf R_\flat}$, $s_1\ge \frac 12$
and $s_2>\frac 12$.
On the other hand, the conditions
on $a$ in Theorems \ref{Thm:GSPseudoDifficultDirection}
and \ref{Thm:GSPseudoEasyDirection} means exactly that
$(z,w)\mapsto e^{(z,w)}a(z,w)$ belongs to the spaces in
\eqref{Eq:SpacesInThms}, depending on the choice between
\eqref{Eq:GSPseudoSymbCond1} and
\eqref{Eq:GSPseudoSymbCond2}, and the condition on $r$.
\end{rem}
%
%

\par

\begin{rem}
Let $s\in \mathbf R_\flat$ be such that $s\le \frac 12$.
By similar arguments as in the proofs of Theorems
\ref{Thm:GSPseudoDifficultDirection}
and \ref{Thm:GSPseudoEasyDirection}, one may also characterize
linear and continuous operators from $\maclA _s'(\cc d)$ to
$\maclA _s(\cc d)$, and from $\maclA _{0,s}'(\cc d)$ to
$\maclA _{0,s}(\cc d)$ as operators of the form
$\op _{\mathfrak V}(a)$ for suitable conditions on $a$. The details are left
for the reader.
\end{rem}

\par

\section{Operators with kernels and symbols in
mixed weighted Lebesgue spaces}\label{sec3}

\par

In this section we focus on operators in the previous section,
whose kernels should belong to $\wideparen A(\cc d \times \cc d)$
and obey certain mixed norm estimates of Lebesgue types. We deduce
continuity properties of such operators when acting between
suitable Lebesgue spaces of analytic functions. (See Theorems
\ref{Thm:LebOpCont}--\ref{Thm:LebOpCont3}.)
Thereafter we show that our results can be used to regain
well-known and sharp continuity results
in \cite{Toft19} for pseudo-differential operators with symbols
in modulation spaces when acting on other modulation spaces.
(See Theorems \ref{Thm:PseudoModCont} and
\ref{Thm:PseudoModCont2}.)
A key step here is to deduce an explicit formula which relates
the short-time Fourier transform
of the symbol to a real pseudo-differential operator $\op (a)$
with the Bargmann transform of the kernel to $\op (a)$.
(See Lemma \ref{Lemma:STFTCompKernelSymbolTransfer}.)

%
%
%
%
%

\par

We shall consider Lebesgue norm conditions of
matrix pull-backs of the involved kernels. Let
\begin{equation}\label{Eq:MatrixCond}
\begin{aligned}
C_{jk}^l &\in \GL (d,\mathbf R),\quad
C_{jk} =
\left (
\begin{matrix}
C_{j,k}^1 & C_{j,k}^2
\\[1ex]
C_{j,k}^3 & C_{j,k}^4
\end{matrix}
\right )
\in \GL (2d,\mathbf R),
\\[1ex]
C &=
\left (
\begin{matrix}
C_{11} & C_{12}
\\[1ex]
C_{21} & C_{22}
\end{matrix}
\right )\in \GL (4d,\mathbf R),\quad
j,k,l\in \mathbf Z_+ ,
\end{aligned}
\end{equation}
and let
\begin{equation}\label{Eq:CondKernels}
\begin{gathered}
U_d(x+i\xi ) = (x,\xi )\in \rr {2d},
\\[1ex]
U_{d,d}(x_1+i\xi _1,x_2+i\xi _2) = (x_1,\xi _1,x_2,\xi _2)\in \rr {4d},
\\[1ex]
K_{\omega} (z,w )
\equiv
e^{-\frac 12 (|z|^2 +|w|^2)} |K (z,w)| \cdot
\omega (\sqrt 2\, \overline z ,\sqrt 2\, w ),
\\[1ex]
G_{K,C,\omega} = K_\omega \circ U_{d,d}^{-1}\circ C \circ U_{d,d},
\\[1ex]
x,x_j ,\xi ,\xi _j \in \rr d,\quad z,w\in \cc d,\ j=1,2.
\end{gathered}
\end{equation}
We will consider continuity of operators from $A^{\mabfp _1}_{E,(\omega _1)}(\cc d)$ to
$A^{\mabfp _2}_{E,(\omega _2)}(\cc d)$, when $G_{K,C,\omega}$ fullfils suitable
$L^{p,q}(\cc {2d})$ estimates, where the weights fullfil
\begin{equation}\label{Eq:CompKernelWeightCond}
\frac {\omega _2(z)}{\omega _1(w)}\lesssim \omega (z,\overline w),\qquad z,w\in \cc d.
\end{equation}
Here and in what follows we let $L^{p,q}(\cc {2d})$ and $L^{p,q}_*(\cc {2d})$ be the sets
of all $G\in L^1_{loc}(\cc {2d})$ such that
$$
\nm G{L^{p,q}(\cc {2d})} \equiv \nm {G\circ U_{d,d}}{L^{p,q}(\rr {4d})}
\quad \text{and}\quad
\nm G{L^{p,q}_*(\cc {2d})} \equiv \nm {G\circ U_{d,d}}{L^{p,q}_*(\rr {4d})},
$$
respectively, are finite.
(See also Remark \ref{Rem:SpaceNotations}.)

\par

The involved Lebesgue exponents should satisfy
\begin{equation}\label{Eq:CondExps}
\begin{gathered}
\frac 1{\mabfp _1} - \frac 1{\mabfp _2} = 1-\frac 1p-\frac 1q ,
\qquad  q\le \mabfp _2\le p,
\\[1ex]
p,q\in [1,\infty ] ,\quad \mabfp _1,\mabfp _2 \in [1,\infty ]^{2d}.
\end{gathered}
\end{equation}
We need that $C$ and $C_{jk}$ above
should satisfy
\begin{align}
\det (C) \det (C_{11}C_{21}) &\neq 0\label{Eq:MatrixCond1}
\intertext{or}
\det (C) \det (C_{12}C_{22}) &\neq 0.\label{Eq:MatrixCond2}
\end{align}
In \eqref{Eq:CondExps} and in what follows we use the convention
$$
\frac 1{\mabfp}
=
\left (
\frac 1{p_1},\dots ,\frac 1{p_d}
\right ),
\quad
p_0\le \mabfp ,
\quad
\mabfp \le p,
\quad
q_0< \mabfq ,
\quad
\mabfq < q
\quad \text{and}\quad
\mabfr =r,
$$
when
$$
\mabfp = (p_1,\dots ,p_d),
\quad
\mabfq = (q_1,\dots ,q_d),
\quad
\mabfr = (r_1,\dots ,r_d)
$$
belong to $[1,\infty ]^d$ and $p,q,r,p_0,q_0,r_0\in [1,\infty ]$
satisfy
$$
p_0\le p_k ,
\quad
p_k \le p,
\quad
q_0< q_k ,
\quad
q_k < q
\quad \text{and}\quad
r_k =r,
\quad k\in \{ 1,\dots ,d \} .
$$
%

\par

\begin{rem}\label{Remark:MatrixCond}
We notice that \eqref{Eq:MatrixCond}--\eqref{Eq:MatrixCond2}
implies that $C$ is invertible and
that one of the following conditions hold true:
\begin{enumerate}
\item both $C_{11}$ and $C_{21}$ are invertible;

\vrum

\item both $C_{12}$ and $C_{22}$ are invertible.
\end{enumerate}
If (1) holds, then
\begin{multline*}
\det
\left (
\begin{matrix}
C_{11}   & C_{12}
\\[1ex]
C_{21}   & C_{22}
\end{matrix}
\right )
\asymp
\det
\left (
\begin{matrix}
I_{2d}   & C_{11}^{-1}C_{12}
\\[1ex]
I_{2d}   & C_{21}^{-1}C_{22}
\end{matrix}
\right )
\\[1ex]
=
\det
\left (
\begin{matrix}
I_{2d}   & C_{11}^{-1}C_{12}
\\[1ex]
0   & C_{21}^{-1}C_{22}-C_{11}^{-1}C_{12}
\end{matrix}
\right )
=
\det (C_{21}^{-1}C_{22}-C_{11}^{-1}C_{12}).
\end{multline*}
Here and in what follows, $I=I_d$ is the $d\times d$ identity matrix.
From these computations it follows that
\begin{alignat*}{2}
&C_{11}^{-1}C_{12}-C_{21}^{-1}C_{22}, & \qquad &C_{11}^{-1}C_{12}-C_{22}C_{21}^{-1}
\\[1ex]
&C_{12}C_{11}^{-1}-C_{21}^{-1}C_{22}, & \qquad &C_{12}C_{11}^{-1}-C_{22}C_{21}^{-1}
\intertext{are invertible when (1) holds, and}
&C_{12}^{-1}C_{11}-C_{22}^{-1}C_{21}, & \qquad &C_{12}^{-1}C_{11}-C_{21}C_{22}^{-1}
\\[1ex]
&C_{11}C_{12}^{-1}-C_{22}^{-1}C_{21}, & \qquad &C_{11}C_{12}^{-1}-C_{21}C_{22}^{-1}
\end{alignat*}
are invertible when (2) holds.
\end{rem}

\par

\begin{rem}
Let $C_0\in \GL (d,\mathbf C)$ and let $U_d$ be the same as in \eqref{Eq:CondKernels}.
Then the matrix $U_d\circ C_0\circ U_d^{-1}$ in $\GL (2d,\mathbf R)$
which corresponds to $C_0$ is given by
\begin{equation}\label{Eq:RealImageC}
\left (
\begin{matrix}
\operatorname{Re} (C_0) & -\operatorname{Im} (C_0)
\\[1ex]
\operatorname{Im} (C_0) & \phantom{-}\operatorname{Re} (C_0)
\end{matrix}
\right ).
\end{equation}
Obviously, the map which takes $C_0\in \GL (d,\mathbf C)$ into the matrix
\eqref{Eq:RealImageC} in $\GL (2d,\mathbf R)$
is injective, but not bijective. In this way we identify $\GL (d,\mathbf C)$
with the set of all matrices in $\GL (2d,\mathbf R)$ which are given by
\eqref{Eq:RealImageC} for some $C_0\in \GL (d,\mathbf C)$.

\par

If $C_{jk}\in \GL (2d,\mathbf R)$ and $T_{jk}=U_d^{-1}\circ C_{jk}\circ U_d$, then
$G_{K,C,\omega}$ in \eqref{Eq:CondKernels} is given by
$$
G_{K,C,\omega} (z,w) = K_{\omega} (T_{11}(z)+T_{12}(w),T_{21}(z)+T_{22}(w)), \quad
z,w \in \cc {d}.
$$

\par

If more restricted, $C_{jk}$ can be identified as matrices in $\GL (d,\mathbf C)$
as above, for $j,k\in \{ 1,2 \}$, then
$G_{K,C,\omega}$ in \eqref{Eq:CondKernels}
is given by
$$
G_{K,C,\omega} (z,w)=K_{\omega} (C_{11}z+C_{12}w,C_{21}z+C_{22}w), \quad
z,w \in \cc {d},
$$
for such choices of $C$.
\end{rem}

\par

\begin{thm}\label{Thm:LebOpCont}
Let $E$ be an ordered basis for $\rr {2d}$,
$\omega _1$ and $\omega _2$ be weights on $\cc d$, $\omega$
be a weight on $\cc d\times \cc d$ such that \eqref{Eq:CompKernelWeightCond} holds,
and let $\mabfp _1$, $\mabfp _2$,
$p$ and $q$ be as in \eqref{Eq:CondExps}. Also let
$C\in \GL (4d,\mathbf R)$ be such that \eqref{Eq:MatrixCond}
holds, $K\in \wideparen A(\cc d\times \cc d)$, and let
$G_{K,C,\omega}$ be as in \eqref{Eq:CondKernels}.
Then the following is true:
\begin{enumerate}
\item if \eqref{Eq:MatrixCond1} holds and
$G_{K,C,\omega}(z,w)\in L^{p,q}(\cc {2d})$,
then $T_K$ in \eqref{Kmap2}
from $\maclA _{\flat _1}(\cc d)$ to $A(\cc d)$ is uniquely extendable to a
continuous mapping from $A^{\mabfp _1}_{E,(\omega _1)}(\cc d)$ to
$A^{\mabfp _2}_{E,(\omega _2)}(\cc d)$,
and
\begin{equation}\label{Eq:OpKernEst1}
\nm {T_KF}{A^{\mabfp _2}_{E,(\omega _2)}}
\lesssim
\nm {G_{K,C,\omega}}{L^{p,q}}\nm {F}{A^{\mabfp _1}_{E,(\omega _1)}},
\quad F\in A^{\mabfp _1}_{E,(\omega _1)}(\cc d)\text ;
\end{equation}

\vrum

\item if \eqref{Eq:MatrixCond2} holds and $G_{K,C,\omega}
\in L^{q,p}_*(\cc {2d})$,
then $T_{K}$ in \eqref{Kmap2}
from $\maclA _{\flat _1}(\cc d)$ to $A(\cc d)$ is uniquely extendable to a
continuous mapping from $A^{\mabfp _1}_{E,(\omega _1)}(\cc d)$ to
$A^{\mabfp _2}_{E,(\omega _2)}(\cc d)$,
and
\begin{equation}\label{Eq:OpKernEst2}
\nm {T_{K}F}{A(\omega _2,L^{\mabfp _2})}
\lesssim
\nm {G_{K,C,\omega}}{L^{q,p}_2}\nm {F}{A^{\mabfp _1}_{E,(\omega _1)}},
\quad F\in A^{\mabfp _1}_{E,(\omega _1)}(\cc d)\text .
\end{equation}
\end{enumerate}
\end{thm}

\par

\begin{proof}
We only prove (1). The assertion (2) follows by similar arguments and is left
for the reader. Let
$$
G_{K,C,\omega ,p}(w) \equiv \nm {G_{K,C,\omega }(\cdo ,w)}{L^p(\cc d)},
\quad w\in \cc d.
$$
Then $\nm {G_{K,C,\omega}}{L^{p,q}(\cc {2d})} \asymp \nm {G_{K,p,\omega}}{L^q(\cc d)}$.
Also let $K_\omega$ be as in \eqref{Eq:CondKernels},
$F\in A^{\mabfp _1}_{E,(\omega _1)}(\cc d)$ and
$H\in B^{\mabfp _2}_{E,(1/\omega _2)}(\cc d)$, and set
\begin{align*}
F_{\omega _1}(w) &\equiv |F(w)|e^{-|w|^2/2}\omega _1(\sqrt 2\, \overline w), \quad w \in \cc {d}
\intertext{and}
H_{\omega _2}(z) &\equiv |H(z)|e^{-|z|^2/2}/\omega _2(\sqrt 2\, \overline z), \quad z \in \cc {d}.
\end{align*}

\par

By H{\"o}lder's inequality we get
\begin{multline*}
|(TF,H)_{B^2}| \asymp
\left |
\int _{\cc d}(TF)(z)\overline{H(z)}e^{-|z|^2}\, d\lambda (z)
\right |
\\[1ex]
\le
{\iint} _{\cc {2d}}K_\omega (z,w)F_{\omega _1}(w)H_{\omega _2}(z)
\, d\lambda (z)d\lambda (w)
\\[1ex]
=
{\iint} _{\cc {2d}}G_{K,C,\omega} (z,w)\Phi (z,w)\, d\lambda (z)d\lambda (w)
\\[1ex]
\lesssim
\nm {G_{K,C,\omega}}{L^{p,q}} \nm {\Phi}{L^{p',q'}},
\end{multline*}
where
\begin{multline*}
\Phi (x+i\xi ,y+i\eta )
\\[1ex]
=
F_{\omega _1}(U_d^{-1}(C_{21}(x,\xi )+C_{22}(y,\eta )))
H_{\omega _2}(U_d^{-1}(C_{11}(x,\xi )+C_{12}(y,\eta ))),
\\[1ex]
x,y,\xi ,\eta \in \rr {d}.
\end{multline*}
Here we identify $(x,\xi )\in \rr {2d}$ by corresponding $2d\times 1$-matrix
$\big ( \, \begin{matrix} x \\[-0.3ex] \xi \end{matrix} \, \big )$,
as usual.

\par

We need to estimate $\nm {\Phi}{L^{p',q'}}$, and start with reformulating
$\nm {\Phi (\cdo ,w)}{L^{p'}}$.  For $\nm {\Phi (\cdo ,w)}{L^{p'}}$ we take
$$
(x,\xi )\mapsto C_{21}((x,\xi )+C_{21}^{-1}C_{22}(y,\eta ))
$$
as new variables of integration, and get
$$
\nm {\Phi (\cdo ,w)}{L^{p'}} \asymp \nm {F_{\omega _1}
\cdot H_{\omega _2}(U_d^{-1}(B_1(\cdo  -B_2(y,\eta )))}{L^{p'}},\quad w=y+i\eta ,
$$
where $B_1$ and $B_2$ are the matrices
$$
B_1= C_{11}C_{21}^{-1}\in \GL (2d,\mathbf R)
\quad \text{and}\quad
B_2 = C_{21}C_{11}^{-1}C_{12}-C_{22} \in \GL (2d,\mathbf R),
$$
which are invertible due to Remark \ref{Remark:MatrixCond} and the assumptions.
Hence, for $F_{\omega _1}^0 = F_{\omega _1}(U_d^{-1} \cdo )$
and $H_{\omega _2}^0 = H_{\omega _2}\circ U_d^{-1} \circ (-B_1)$ we have
\begin{equation}\label{Eq:PhiH0Reform}
\nm {\Phi (\cdo ,w)}{L^{p'}}
\asymp
\left (\left ( |F_{\omega _1}^0| ^{p'}*
|H_{\omega _2}^0|^{p'} \right ) (B_2 (y,\eta )) \right )^{\frac 1{p'}},
\quad w=y+i\eta.
\end{equation}

\par

If $\mabfr _1=\mabfp _1/p'$ and $\mabfr _2=\mabfp _2'/p'$,
then it follows from \eqref{Eq:CondExps} that
$$
\frac 1{\mabfr _1}+\frac 1{\mabfr _2}= 1+\frac {p'}{q'},
\quad \text{and}\quad
\mabfr _1,\mabfr _2, \frac {q'}{p'}\ge 1.
$$
Hence, by \eqref{Eq:PhiH0Reform}, that $q'/p'\ge 1$, the fact that $B_2$
is invertible, and H{\"o}lder's and Young's inequalities we obtain
\begin{multline*}
\nm {\Phi}{L^{p',q'}}
\lesssim
\Nm {\left (\left ( |F_{\omega _1}^0| ^{p'}*
|H_{\omega _2}^0|^{p'} \right ) (B_2 \cdo ) \right )^{\frac 1{p'}}}{L^{q'}}
\\[1ex]
\asymp \left (
\Nm {|F_{\omega _1}^0| ^{p'}*
|H_{\omega _2}^0|^{p'} }{L^{q'/p'}_E}
\right )^{\frac 1{p'}}
\\[1ex]
\le
\left (
\nm {|F_{\omega _1}^0| ^{p'}}{L^{\mabfr _1}_E}
\nm {|H_{\omega _2}^0| ^{p'}}{L^{\mabfr _2}_E}
\right )^{\frac 1{p'}}
\asymp
\nm {F_{\omega _1}}{L^{\mabfp _1}_E}\nm {H_{\omega _2}}{L^{\mabfp _2'}_E}
\end{multline*}
and the right-hand side of \eqref{Eq:OpKernEst1} follows by taking the supremum
over all such $H$ with $\nm H{B^{\mabfp _2'}_{E,(1/\omega _2)}}\le 1$.

\par

The existence of extension now follows from Hahn-Banach's theorem. By these
estimates it also follows that
$$
(z,w)\mapsto
K_{\omega}(z,w)F_{\omega _1}(w)\overline{H_{\omega _2}(z)}
$$
belongs to $L^1(\cc d\times \cc d)$, and the uniqueness is a straight-forward application
of Lebesgue's theorem.
\end{proof}

\par

For corresponding pseudo-differential operator with symbol $a$, the
kernel is given by $K(z,w)=e^{(z,w)}a(z,w)$. By straight-forward
computations it follows that $G_{K,C,\omega ,p}$ takes the form
\begin{align}
G_{K,C,\omega ,p}(w) &= \nm {a_\omega (\cdo ,w)}{L^p(\cc d)},\quad w\in \cc d,
\notag
\intertext{when $C_{11}=C_{12}=C_{22}=I_{2d}$ and $C_{21}=0$, or
$C_{11}=C_{12}=C_{21}=I_{2d}$ and $C_{22}=0$, where}
a_\omega (z,w) &= e^{-\frac 12 \cdot |z|^2}a(z+w,w)
\omega (\sqrt 2\, (\overline{z+w}), \sqrt 2\, w)
\label{Eq:aOmegaDef1}
\intertext{or}
a_\omega (z,w) &= e^{-\frac 12 \cdot |w|^2}a(z+w,z)
\omega (\sqrt 2\, (\overline{z+w}), \sqrt 2\, z). \label{Eq:aOmegaDef2}
\end{align}
Hence, Theorem \ref{Thm:LebOpCont} gives the following.

\par

\begin{thm}\label{Thm:LebOpCont2}
Let $\omega _1$ and $\omega _2$ be weights on $\cc d$, $\omega$
be a weight on $\cc d\times \cc d$ such that \eqref{Eq:CompKernelWeightCond} holds,
$\mabfp _1$, $\mabfp _2$, $p$ and $q$ be as in \eqref{Eq:CondExps}.
Also let $a\in \wideparen A(\cc d\times \cc d)$ and let $a_\omega$ be
given by \eqref{Eq:aOmegaDef1} or by \eqref{Eq:aOmegaDef2} for $z,w\in \cc d$. If
$a_\omega \in L^{p,q}(\cc {2d})$, then
the operator $\op _{\mathfrak V}(a)$ in \eqref{APsiDO}
from $\maclA _{\flat _1}(\cc d)$ to $A(\cc d)$ is uniquely extendable to a continuous
mapping from $A^{\mabfp _1}_{E,(\omega _1)}(\cc d)$ to
$A^{\mabfp _2}_{E,(\omega _2)}(\cc d)$.
\end{thm}

\par

We also have the following result related to Theorem \ref{Thm:LebOpCont}. Here the
matrix $C$ is given by \eqref{Eq:MatrixCond} with
\begin{equation}\label{Eq:MatrixCondSpec1}
C_{11} = C_{21} =
I_{2d},
\quad
C_{12} =
\left (
\begin{matrix}
0 & 0
\\[1ex]
0 & I_d
\end{matrix}
\right )
\quad \text{and}\quad
C_{22} =
\left (
\begin{matrix}
I_d & 0
\\[1ex]
0 & 0
\end{matrix}
\right )
\end{equation}
which obviously satisfies \eqref{Eq:MatrixCond1}. Also again recall Remark
\ref{Rem:SpaceNotations} for notations.

\par

\begin{thm}\label{Thm:LebOpCont3}
Let $C$ be given by \eqref{Eq:MatrixCond} with $C_{jk}$ given by \eqref{Eq:MatrixCondSpec1},
$\omega _1$ and $\omega _2$ be weights on $\cc d$, $\omega$
be a weight on $\cc d\times \cc d$ such that \eqref{Eq:CompKernelWeightCond} holds,
and let $p,q\in [1,\infty ]$. Also let
$K\in \wideparen A(\cc d\times \cc d)$ and $G_{K,C,\omega}$ be as in
\eqref{Eq:CondKernels}.
If $G_{K,C,\omega}\in L^{p,q}_*(\cc {2d})$, then
then $T_K$ in \eqref{Kmap2}
from $\maclA _{\flat _1}(\cc d)$ to $A(\cc d)$ is uniquely extendable to a
continuous mapping from $A^{p',q'}_{(\omega _1)}(\cc d)$ to
$A^{q,p}_{*,(\omega _2)}(\cc d)$, and
\begin{equation}\label{Eq:OpKernEst1A}
\nm {T_KF}{A^{q,p}_{*,(\omega _2)}}
\lesssim
\nm {G_{K,C,\omega}}{L^{p,q} _*}\nm {F}{A^{p',q'}_{(\omega _1)}},
\quad F\in A^{p',q'}_{(\omega _1)}(\cc d).
\end{equation}
\end{thm}

\par

\begin{proof}
Let $F_{\omega _1}$, $F_{\omega _1}^0$, $H_{\omega _2}$ and $H_{\omega _2}^0$
be the same as in the proof of Theorem \ref{Thm:LebOpCont}, and let
$K_\omega$ be as in \eqref{Eq:CondKernels}. Then
\begin{multline*}
|(T_KF,H)_{B^2}|
\le
\iiiint K_\omega (x,\xi ,y,\eta)
F_{\omega _1}^0(y,\eta ){H}_{\omega _2}^0(x,\xi )\, dxd\xi dyd\eta
\\[1ex]
=
\iiiint G_{K,C,\omega } (x,\xi ,y,\eta ))F_{\omega _1}^0(x+y,\xi )
{H}_{\omega _2}^0(x,\xi +\eta )
\, dxd\xi dyd\eta
\\[1ex]
\le
\nm {G_{K,C,\omega}}{L^{p,q}_*} \nm {\Phi _0}{L^{p'}},
\end{multline*}
where
\begin{multline*}
\Phi _0(x,\xi )
=
\left (
\iint |F_{\omega _1}^0(x+y,\xi )H_{\omega _2}^0(x,\xi +\eta )|^{q'}\, dyd\eta
\right )^{\frac 1{q'}}
\\[1ex]
=
\nm {F_{\omega _1}^0(\cdo ,\xi )}{L^{q'}}\nm {H_{\omega _2}^0(x,\cdo )}{L^{q'}}, \quad
x,\xi \in \rr {d}.
\end{multline*}
Hence,
$$
\nm {\Phi _0}{L^{p'}} = \nm {F_{\omega _1}}{L^{q',p'}_1}\nm {H_{\omega _2}}{L^{p',q'}_2}.
$$
The continuity assertion now follows from these estimates, an application of Hahn-Banach's
theorem and Lebesgue's theorem (cf. the end of the proof of Theorem \ref{Thm:LebOpCont}).
\end{proof}

\par

\begin{rem}
In Theorem \ref{Thm:LebOpCont3}, the matrix $C$ is chosen only as
\eqref{Eq:MatrixCond} and \eqref{Eq:MatrixCondSpec1}, while
Theorem \ref{Thm:LebOpCont} is valid for a whole family of matrices
with the only restriction \eqref{Eq:MatrixCond1} or \eqref{Eq:MatrixCond2}.
On the other hand, by similar arguments, it follows that the conclusions
in Theorem \ref{Thm:LebOpCont3} are still true when, more generally,
$C \in \GL (4d, \mathbf R)$ is of the form
\begin{alignat*}{2}
&\left (
\begin{matrix}
\pm I_d & 0 & \pm I_d & 0
\\[1ex]
0 & \pm I_d & 0 & 0
\\[1ex]
\pm I_d & 0 & 0 & 0
\\[1ex]
0 & \pm I_d & 0 & \pm I_d
\end{matrix}
\right ), &
\qquad
&\left (
\begin{matrix}
\pm I_d & 0 & 0 & 0
\\[1ex]
0 & \pm I_d & 0 & \pm I_d
\\[1ex]
\pm I_d & 0 & \pm I_d & 0
\\[1ex]
0 & \pm I_d & 0 & 0
\end{matrix}
\right ),
\\[1ex]
&\left (
\begin{matrix}
\pm I_d & 0 & \pm I_d & 0
\\[1ex]
0 & 0 & 0 & \pm I_d
\\[1ex]
0 & 0 & \pm I_d & 0
\\[1ex]
0 & \pm I_d & 0 & \pm I_d
\end{matrix}
\right ) &
\quad \text{or}\quad
&\left (
\begin{matrix}
0 & 0 & \pm I_d & 0
\\[1ex]
0 & \pm I_d & 0 & \pm I_d
\\[1ex]
\pm I_d & 0 & \pm I_d & 0
\\[1ex]
0 & 0 & 0 & \pm I_d.
\end{matrix}
\right ) ,
\end{alignat*}
for any choice of $\pm$ at each place, provided the mixed Lebesgue conditions on
$G_{K,C,\omega}$ are slightly modified.
\end{rem}

\medspace

In order to apply Theorem \ref{Thm:LebOpCont} to \emph{real}
pseudo-differential operators we have the following.

\par

\begin{lemma}\label{Lemma:STFTCompKernelSymbolTransfer}
Let $\phi (x,\xi )=\pi ^{-\frac d2}e^{i\scal x\xi}e^{-\frac 12( |x |^2+|\xi |^2)}$,
$x,\xi \in \rr d$,
$a\in \maclH _{\flat _1}'(\rr {2d})$ and let $K_a$ be the kernel of
$\op (a)$. Then
\begin{multline*}
e^{-\frac 12(|z|^2+|w|^2)}\mathfrak V_{\Theta ,d}K_a(z,w)
\\[1ex]
=
2^{\frac d2}e^{-i(\scal x{\xi -2\eta} +\scal y\eta)}
(V_\phi a)(\sqrt 2 x,-\sqrt 2\eta ,\sqrt 2(\eta -\xi ),\sqrt 2(y-x))
\end{multline*}
when $z=x+i\xi \in \cc d$ and $w=y+i\eta \in \cc d$.
\end{lemma}

\par

\begin{proof}
Let $\phi _0(x,\xi )=\pi ^{-\frac d2}e^{-\frac 12(|x |^2+|\xi |^2)} = e^{-i\scal x\xi}\phi (x,\xi )$.
By formal computations and Fourier's inversion formula we get
\begin{multline*}
(4\pi ^3)^{\frac d2}e^{\frac i2\cdot (\scal x\xi +\scal y\eta)}e^{-\frac 14(|z|^2+|w|^2)}
\mathfrak V_{\Theta ,d}K_a(z/{\sqrt 2},w/{\sqrt 2})
\\[1ex]
=
\iiint a(x_1,\xi _1)e^{i\scal {x_1-y_1}{\xi _1}} \phi _0(x_1-x,y_1-y
)e^{-i(\scal {y_1}\eta - \scal {x_1}\xi )}\, dx_1dy_1d\xi _1
\\[1ex]
=
(2\pi )^{\frac d2}e^{-i\scal y\eta}
\iint a(x_1,\xi _1)\phi _0(x_1-x,\xi _1+\eta )
e^{i\scal {x_1}{\xi _1}}e^{-i(\scal y{\xi _1}-\scal {x_1}\xi)}\, dx_1d\xi _1
\\[1ex]
=
(2\pi )^{\frac d2} e^{i\scal {x-y}\eta}
\iint a(x_1,\xi _1)\phi (x_1-x,\xi _1+\eta )
e^{-i(\scal {x_1}{\eta -\xi} + \scal {y-x}{\xi _1})}\, dx_1d\xi _1
\\[1ex]
=
(4\pi )^{\frac d2} e^{i\scal {x-y}\eta}
(V_\phi a)(x,-\eta ,\eta -\xi ,y-x). 
\qedhere
\end{multline*}
\end{proof}

\par

We can now use the previous lemma and theorems
to obtain mapping properties for pseudo-differential operators with
symbols in modulation spaces. For example, we may combine
Lemma \ref{Lemma:STFTCompKernelSymbolTransfer} and
Theorem \ref{Thm:LebOpCont} to deduce the following result,
which is the same as \cite[Theorem 2.2]{Toft22}. Hence our kernel results
on the Bargmann transform side can be used to regain classical
mapping properties pseudo-differential operators when acting on
modulation spaces.

\par

\begin{thm}\label{Thm:PseudoModCont}
Let $E$ be an ordered basis of $\rr {2d}$, $A\in \GL (\mathbf R,d)$,
$p,q\in [1,\infty ]$ and $\mabfp _1,\mabfp _2\in [1,\infty ]^{2d}$
be as in \eqref{Eq:CondExps}, $\omega _0\in \mascP _E(\rr {4d})$
and $\omega _1,\omega _2 \in \mascP _E(\rr {2d})$ be such that
\begin{equation}\label{Eq:PseudWeightCond}
\frac {\omega _2(x-Ay,\xi +(I-A^*)\eta )}{\omega _1(x+(I-A)y,\xi -A^*\eta )}
\lesssim
\omega _0(x,\xi ,\eta ,y),\quad x,y,\xi ,\eta \in \rr d,
\end{equation}
and let $a\in M^{p,q}_{(\omega _0)}(\rr {2d})$.
Then $\op _A(a)$ from $\Sigma _1(\rr d)$ to $\Sigma _1'(\rr d)$ extends
uniquely to continuous operator from $M^{\mabfp _1}_{E,(\omega _1)}(\rr d)$
to $M^{\mabfp _2}_{E,(\omega _2)}(\rr d)$, and
$$
\nm {\op _A(a)}{M^{\mabfp _1}_{E,(\omega _1)} \to M^{\mabfp _2}_{E,(\omega _2)}}
\lesssim
\nm a{M^{p,q}_{(\omega _0)}}.
$$
\end{thm}

\par

\begin{proof}
By \eqref{calculitransform} and Proposition \ref{Prop:ExpOpSTFT} we may
assume that $A=0$.
Let
$$
C_{11}=C_{21}=C_{22}=I_{2d},\quad C_{12}=0
$$
and let $\omega$ be given by
$$
\omega (x,\xi ,y,\eta ) = \omega _0(x,-\eta ,\xi +\eta ,y-x),\quad x,y,\xi ,\eta \in \rr d,
$$
which we identify with
$$
\omega (z,w),\qquad z=x+i\xi \in \cc {d}, w=y+i\eta \in \cc {d}.
$$
Then it follows by straight-forward computations that
\eqref{Eq:PseudWeightCond} is the same as
\eqref{Eq:CompKernelWeightCond}.
Furthermore, let $K_0=\mathfrak V_{\Theta ,d}K_a$, where $K_a$ is the kernel
of the operator $\op (a)$. Then it follows from Lemma
\ref{Lemma:STFTCompKernelSymbolTransfer}
and straight-forward computations that if
$$
H_{a,\omega _0}(x,\xi ,\eta ,y) = |V_\phi a (x,\xi ,\eta ,y)|\cdot
\omega _0 (x,\xi ,\eta ,y),,\quad x,y,\xi ,\eta \in \rr d,
$$
then
\begin{equation}\label{EqHaGKRel}
\begin{gathered}
H_{a,\omega _0}(\sqrt 2\, x , -\sqrt 2 (\xi +\eta ), \sqrt 2 \, \eta ,\sqrt 2 \, y)
\asymp
G_{K_0,C,\omega }(z,w),
\\[1ex]
z=x+i\xi \in \cc {d},\ w=y+i\eta \in \cc {d}.
\end{gathered}
\end{equation}

\par

By first applying the $L^p$-norm on \eqref{EqHaGKRel} with respect
to $x$ and $\xi$, and thereafter applying the $L^q$-norm with respect to $y$
and $\eta$, we get
$$
\nm {G_{K_0,C,\omega }}{L^{p,q}(\cc {2d})}
\asymp
\nm {H_{a,\omega _0}}{L^{p,q}(\rr {4d})}
\asymp
\nm a{M^{p,q}_{(\omega _0)}}<\infty .
$$

\par

Hence, the assumptions in Theorem \ref{Thm:LebOpCont} are fullfiled,
and we conclude that the operator $T_{K_0}$ with kernel $K_0$ is
continuous from $A^{\mabfp _1}_{E,(\omega _1)}(\cc d)$ to
$A^{\mabfp _2}_{E,(\omega _2)}(\cc d)$. The asserted continuity for $\op (a)$
is now a consequence of the commutative diagram
\begin{equation}  \label{Eq:CommDiagram2}
\begin{CD}
M^{\mabfp _1}_{E,(\omega _1)}(\rr d)    @> \op (a) >> M^{\mabfp _2}_{E,(\omega _2)}(\rr d)
\\
@V {{ \mathfrak V}_d} VV        @VV{{ \mathfrak V}_d}V
\\
A^{\mabfp _1}_{E,(\omega _1)}(\cc d)     @>>T_{K_0}>
A^{\mabfp _2}_{E,(\omega _2)}(\cc d). \qedhere
\end{CD}
\end{equation}
\end{proof}

\par

The next result extends \cite[Theorem 3.3]{To10} and
follows by similar arguments as in the previous proof, using
Theorem \ref{Thm:LebOpCont3} instead of Theorem \ref{Thm:LebOpCont}.
The details are left for the reader.

\par

\begin{thm}\label{Thm:PseudoModCont2}
Let $\omega _1$ and $\omega _2 \in \mascP _E(\rr {2d})$, $\omega \in
\mascP _E(\rr {4d})$ be such that
$$
\frac {\omega _2(x,\xi +\eta )}{\omega _1(x+y,\xi )} \lesssim \omega _0(x,\xi ,\eta ,y),
,\quad x,y,\xi ,\eta \in \rr d,
$$
let $p,q\in [1,\infty ]$, and let $a\in W^{p,q}_{(\omega _0)}(\rr {2d})$.
Then $\op _0(a)$ from $\Sigma _1(\rr d)$ to $\Sigma _1'(\rr d)$ is uniquely
extendable to a continuous mapping from $M_{(\omega _1)}^{q',p'}(\rr d))$ to
$W_{(\omega _2)}^{p,q}(\rr d))$, and
$$
\nm {\op (a)f}{W_{(\omega _2)}^{p,q}}
\lesssim
\nm {a}{W^{p,q}_{(\omega _0)}}\nm {f}{M_{(\omega _1)}^{q',p'}},
\quad f\in M_{(\omega _1)}^{q',p'}(\rr d)).
$$
\end{thm}

\par

\begin{rem}\label{Rem:Some Extensions}
Let $E$, $p$, $q$, $\mabfp _j$, $\phi _0$ and $\phi$ be the same as in
Lemma \ref{Lemma:STFTCompKernelSymbolTransfer},
Theorem \ref{Thm:PseudoModCont} and their proofs. Also let $A\in \GL (\mathbf R,d)$
and let $\phi _A = e^{i\scal {AD_\xi }{D_x}}\phi $. Then the condition on $a$
in Theorem \ref{Thm:PseudoModCont} is that $\nm {V_{\phi _0}a\cdot
\omega _0}{L^{p,q}}<\infty$. In view of \cite{CaTo,CarWal,Tr} and
Proposition \ref{Prop:BasicModSpaceProperties} (2), the previous condition
is the same as $\nm {V_{\phi _A}a\cdot \omega _0}{L^{p,q}}<\infty$ because
$\omega _0$ is moderate.

\par

We observe that all weights in
Theorem \ref{Thm:PseudoModCont} are moderate, while
there are no such assumptions or other restrictions on the involved
weight functions in Theorem \ref{Thm:LebOpCont}. Since the latter result
is used to prove the former one, a natural question is wether
Theorem \ref{Thm:PseudoModCont} can be extended to broader classes of
weight functions. In view of Remark \ref{Rem:BasicModSpaceProperties},
it is evident that the imposing moderate conditions on weights might
in some context be considered as strong restrictions.

\par

The answer on this question is affirmative in the sense that for suitable
modifications, the moderate conditions on the weights in
Theorem \ref{Thm:PseudoModCont} can be removed.

\par

In fact, let $\maclH _{\flat _1}^A(\rr {2d})$ be the modification of
$\maclH _{\flat _1}(\rr {2d})$, given by
$$
\maclH _{\flat _1}^A(\rr {2d})
=
\sets {e^{i\scal {AD_\xi }{D_x}}(e^{i\scal x\xi }a)}{a\in \maclH _{\flat _1}(\rr {2d})},
$$
$(\maclH _{\flat _1}^A)'(\rr {2d})$ be the dual of $\maclH _{\flat _1}^A(\rr {2d})$,
$\omega _1,\omega _2$ be weights on $\rr {2d}$ and let $\omega _0$
be a weight on $\rr {4d}$ such that \eqref{Eq:PseudWeightCond} holds.
Then it follows from the proof of Theorem \ref{Thm:PseudoModCont} that the following
is true:
\begin{itemize}
\item
if $a\in (\maclH _{\flat _1}^A)'(\rr {2d})$, then $V_{\phi _A}a$ makes sense as
a smooth function;

\vrum

\item
if $a\in (\maclH _{\flat _1}^A)'(\rr {2d})$ satisfies
$\nm {V_{\phi _A}a\cdot \omega _0}{L^{p,q}}<\infty$, then
$\op _A(a)$ from $\maclH _{\flat _1}(\rr d)$ to $\maclH _{\flat _1}'(\rr d)$ extends
uniquely to a continuous operator from $M^{\mabfp _1}_{E,(\omega _1)}(\rr d)$
to $M^{\mabfp _2}_{E,(\omega _2)}(\rr d)$.
\end{itemize}

\par

In similar ways, Theorem \ref{Thm:PseudoModCont2} can be extended to
permit more general weight classes.
\end{rem}

\par


\begin{thebibliography}{99}
\bibitem{B1} V. Bargmann \emph{On a Hilbert space of analytic
functions and an associated integral transform}, Comm. Pure
Appl. Math., \textbf{14} (1961), 187--214.

\bibitem{B2} V. Bargmann \emph{On a Hilbert space of analytic functions
 and an associated integral transform. Part II. A family of related
function spaces. Application to distribution theory.}, Comm. Pure
Appl. Math., \textbf{20} (1967), 1--101.

\bibitem{Bau} W. Bauer \emph{Berezin-Toeplitz quantization and
composition formulas}, {J. Funct. Anal.}, \textbf{256}
(2007), 3107--3142.

\bibitem{Berezin71} F. A. Berezin
\emph{Wick and anti-{W}ick symbols of operators},
Mat. Sb. (N.S.), \textbf{86} (1971), 578--610.


\bibitem{CaTo}
M. Cappiello, J. Toft, \emph{Pseudo-differential operators
in a Gelfand--Shilov setting},  Math. Nachr. \textbf{290} (2017), 738--755.

\bibitem{CarWal}
E. Carypis, P. Wahlberg, \emph{Propagation of exponential phase
space singularities for Schr{\"o}dinger equations with quadratic Hamiltonians},
J. Fourier Anal. Appl. \textbf{23} (2017), 530--571.

\bibitem{CheSigTo2} Y. Chen, M. Signahl, J. Toft
\emph{Factorizations and singular value estimates of
operators with Gelfand-Shilov and Pilipovi{\'c} kernels},
J. Fourier Anal. Appl. \textbf{24} (2018), 666--698.

\bibitem{CPRT10} E. Cordero, S. Pilipovi\'c, L. Rodino, N. Teofanov
\emph{Quasianalytic Gelfand-Shilov spaces with applications
to localization operators}, Rocky Mt. J. Math. \textbf{40} (2010),
1123-1147.

\bibitem{Fe0} H. G. Feichtinger \emph{Banach spaces of
distributions of Wiener's type and interpolation, \rm {in:
Ed. P. Butzer, B. Sz. Nagy and E. G{\"o}rlich (Eds), Proc.
Conf. Oberwolfach, Functional Analysis and Approximation,
August 1980}}, Int. Ser. Num. Math. \textbf{69} Birkh{\"a}user Verlag,
Basel, Boston, Stuttgart, 1981, pp. 153--165.

\bibitem{F1}  H. G. Feichtinger \emph{Modulation spaces on locally
compact abelian groups. Technical report}, {University of
Vienna}, Vienna, 1983; also in: M. Krishna, R. Radha,
S. Thangavelu (Eds) Wavelets and their applications, Allied
Publishers Private Limited, NewDehli Mumbai Kolkata Chennai Hagpur
Ahmedabad Bangalore Hyderbad Lucknow, 2003, pp. 99--140.

\bibitem{FeGaTo2} C. Fernandez, A. Galbis, J. Toft
\emph{The Bargmann transform and powers of harmonic oscillator on
Gelfand-Shilov subspaces}, RACSAM \textbf{111} (2017), 1--13.


%

\bibitem{Gc2} {K. Gr{\"o}chenig} \newblock \emph{Foundations of
Time-Frequency Analysis},
\newblock Birkh{\"a}user, Boston, 2001.

\bibitem{Groch} K. Gr{\"o}chenig \emph{Weight functions in time-frequency analysis
\rm {in: L. Rodino, M. W. Wong (Eds) Pseudodifferential
Operators: Partial Differential Equations and Time-Frequency Analysis}},
Fields Institute Comm., \textbf{52} 2007, pp. 343--366.

\bibitem{GZ} K. Gr{\"o}chenig, G. Zimmermann
\emph{Spaces of test functions via the STFT} J. Funct. Spaces Appl. \textbf{2}
(2004), 25--53.

\bibitem{Ho1} L. H{\"o}rmander \emph{The Analysis of Linear
Partial Differential Operators}, vol {I--III},
Springer-Verlag, Berlin Heidelberg NewYork Tokyo, 1983, 1985.

\bibitem{LozPer} Z. Lozanov Crvenkovi{\'c}, D. Peri{\v s}i{\'c}
\emph{Hermite expansions of elements of Gelfand Shilov spaces
in quasianalytic and non quasianalytic case},
Novi Sad J. Math. \textbf{37} (2007), 129--147.

\bibitem{Pil1} S. Pilipovi\'c \emph{Generalization of Zemanian spaces
of generalized functions which
have orthonormal series expansions},
SIAM J. Math. Anal. \textbf{17} (1986), 477--484.

\bibitem{Pil2}
S. Pilipovi\'c \emph{Tempered ultradistributions},
Boll. U.M.I. \textbf{7} (1988), 235--251.

\bibitem{RS} M. Reed and B. Simon  \emph{Methods of modern
mathematical physics}, Academic Press, London New York, 1979.

\bibitem{RuTo} M. Ruzhansky, N. Tokmagambetov
\emph{Nonharmonic analysis of boundary value problems},
Int. Math. Res. Notices \textbf{12} (2016), 3548--3615.

\bibitem{Teof}
N. Teofanov
\emph{Ultradistributions and time-frequency analysis
{\rm {in: P. Boggiatto, L. Rodino, J. Toft, M. W. Wong (eds)}}
Pseudo-differential operators and related topics},
Operator Theory: Advances and Applications \textbf{164},
Birkh{\"a}user, Basel, 2006, pp. 173--192.

\bibitem{Teof2}
N. Teofanov,
\emph{Gelfand-Shilov spaces and localization operators},
Funct. Anal. Approx. Comput. \textbf{7} (2015), 135--158.


\bibitem{To8} J. Toft \emph{Continuity
properties for modulation spaces with applications to
pseudo-differential calculus, II}, {Ann. Global Anal. Geom.}
\textbf{26} (2004), 73--106.

\bibitem{To10} J. Toft \emph{Pseudo-differential operators with symbols in
modulation spaces}, in: B.-W. Schulze, M. W. Wong (Eds),
Pseudo-Differential Operators: Complex Analysis and Partial
Differential Equations, Operator Theory  Advances and Applications
\textbf{205}, Birkh{\"a}user Verlag, Basel, 2010, pp. 223--234.

\bibitem{To11} J. Toft \emph{The Bargmann transform on modulation and
Gelfand-Shilov spaces, with applications to Toeplitz and pseudo-differential
operators},  J. Pseudo-Differ. Oper. Appl.  \textbf{3}  (2012), 145--227.


\bibitem{Toft18} J. Toft
\emph{Images of function and distribution spaces under the
Bargmann transform}, J. Pseudo-Differ. Oper. Appl. \textbf{8} (2017), 83--139.

\bibitem{Toft19} J. Toft \emph{Continuity and compactness for pseudo-differential operators
with symbols in quasi-Banach spaces or H{\"o}rmander classes}, Anal. Appl.
\textbf{15} (2017), 353--389.

\bibitem{Toft22} J. Toft \emph{Matrix parameterized pseudo-differential
calculi on modulation spaces {\rm {in: M. Oberguggenberger, J. Toft,
J. Vindas, P. Wahlberg (eds)}} Generalized
Functions and Fourier Analysis}, Operator Theory: Advances and
Applications \textbf{260}, Birkh{\"a}user, Basel Heidelberg NewYork
Dordrecht London, pp. 215--235.


\bibitem{Tr} G. Tranquilli \emph{Global normal forms and global properties in
function spaces for second order Shubin type operators} PhD Thesis, 2013.
\end{thebibliography}
\end{document}